\documentclass[12pt]{amsart}
 \usepackage[dvips]{graphicx}
\usepackage[all]{xypic}
\usepackage{amsmath,graphics}
\usepackage{amsfonts,amssymb}
\theoremstyle{plain}
\newtheorem*{theorem*}{Theorem}
\newtheorem*{lemma*} {Lemma}
\newtheorem*{corollary*} {Corollary}
\newtheorem*{proposition*} {Proposition}
\newtheorem*{conjecture*}{Conjecture}

\newtheorem{theorem}{Theorem}[section]
\newtheorem{lemma}[theorem]{Lemma}
\newtheorem{corollary}[theorem]{Corollary}
\newtheorem{proposition}[theorem]{Proposition}

\newcommand{\dimm}{\operatorname{dim}}

\theoremstyle{remark}
\newtheorem*{remark}{Remark}

\newtheorem*{claim}{Claim}

\theoremstyle{definition}

\def\eps{\epsilon}

\def\s{\sigma}

\def\gl{\mbox{GL}}
\def\Q{\Bbb{Q}}
\def\F{\Bbb{F}}
\def\id{\mbox{id}}

\def\Z{\Bbb{Z}}
\def\R{\Bbb{R}}

\def\genus{\mbox{genus}}

\def\part{\partial}

\def\a{\alpha}

\def\g{\gamma}

\def\degtaum{\deg(\tau(M,\phi,\a))}
\def\n{\mathfrak{n}}
\def\m{\mathfrak{m}}
\def\tor{\mbox{Tor}}
\def\bp{\begin{pmatrix}}

\def\sm{\setminus}
\def\ep{\end{pmatrix}}
\def\bn{\begin{enumerate}}

\def\Hom{\mbox{Hom}}

\def\en{\end{enumerate}}
\def\ba{\begin{array}}
\def\ea{\end{array}}
\def\be{\begin{equation}}
\def\ee{\end{equation}}

\def\s{\sigma}
\def\a{\alpha}
\def\ol{\overline}
\def\b{\beta}
\def\ti{\tilde}

\def\fr12{\frac{1}{2}}

\def\i{\iota}
\def\at{\overline{\a}}

\def\im{\mbox{Im}}

\def\ker{\mbox{Ker}}

\def\hom{\mbox{Hom}}

\def\v{\varphi}

\def\ext{\mbox{Ext}}

\def\H{\mathcal{H}}

\def\dim{\mbox{dim}}

\def\tpm {[t^{\pm 1}]}

 \def\twialexi{\Delta_{i}^{\a}(t)}
 
 \def\titwialexi{\ti{\Delta}_{i}^{\a}(t)}
 \def\degtwialexzero{\deg\left(\Delta_{0}^{\a}(t)\right)}

\def\degtwialexi{\deg\left(\Delta_{i}^{\a}(t)\right)}
\def\twialexzero{\Delta_{0}^{\a}(t)} 
\def\twialexg0{\Delta^{G}_{0}(t)}

\def\degtwialexzero{\deg\left(\Delta_{0}^{\a}(t)\right)}

\def\degtwialexg0{\deg\left(\Delta^{G}_{0}(t)\right)} 
 \def\twialexgk0{\Delta^{G}_{0}(t)}

\def\twialextwo{\Delta_{2}^{\a}(t)}

\def\degtwialextwo{\deg\left(\Delta_{2}^{\a}(t)\right)} \def\twialex{\Delta_{1}^{\a}(t)}
 \def\twialexg{\Delta^{G}_{1}(t)}
\def\degtwialex{\deg\left(\Delta_{1}^{\a}(t)\right)}
\def\degtwialexk{\deg\left(\Delta_{1}^{\a}(t)\right)}
\def\degtwialexg{\deg\left(\Delta^{G}_{1}(t)\right)} 
 \def\twialexgk{\Delta^{G}_1(t)}

 \def\alex{\Delta_{1}(t)} \def\alexk{\Delta_{K}(t)}

 \def\twihomi{H_i^{\a}(M;\fkt)}

 \def\twihomphi{H_1^{\a\otimes \phi}(M;\fkt)}
 
 \def\twihomphizero{H_0^{\a\otimes \phi}(M;\fkt)}
 \def\twihomphitwo{H_2^{\a\otimes \phi}(M;\fkt)}

\def\fpt{\F_p[t^{\pm 1}]}

\def\fkt{\F^k[t^{\pm 1}]}
\def\ft{\F[t^{\pm 1}]}
\def\f{\F}
\def\GL{\mbox{GL}}
\def\glfk{\GL(\F,k)}

\def\fk{\F^k}

\def\qt{\Q[t^{\pm 1}]}

\def\K{\Bbb{K}}

\def\tnphi{||\phi||_T}

\def\gen{\mbox{genus}}

\def\oneplusbm{\big(1+b_3(M)\big)}

\begin{document}

\title[Thurston norm, fibered manifolds and twisted Alexander polynomials]
{The Thurston norm, fibered manifolds and twisted Alexander polynomials} \author{Stefan Friedl and Taehee Kim}
\date{\today} \address{Rice University, Houston, Texas, 77005-1892}
\email{friedl@math.rice.edu}\address{ Department of Mathematics,
Konkuk University, Hwayang-dong, Gwangjin-gu, Seoul 143-701,
Korea} \email{tkim@konkuk.ac.kr}
\def\subjclassname{\textup{2000} Mathematics Subject Classification}
\expandafter\let\csname
subjclassname@1991\endcsname=\subjclassname
\expandafter\let\csname
subjclassname@2000\endcsname=\subjclassname \subjclass{Primary
57M27; Secondary 57N10} \keywords{Thurston norm, Twisted Alexander
polynomials, 3-manifolds, Knot genus, Fibered knots}

\begin{abstract}
Every element in the first cohomology group of a 3--manifold is dual to embedded surfaces. The
Thurston norm measures the minimal `complexity' of such surfaces. For instance the Thurston norm
of a knot complement determines the genus of the knot in the 3--sphere. We show that the degrees
of twisted Alexander polynomials give lower bounds on the Thurston norm, generalizing work of
McMullen and Turaev. Our bounds attain their most concise form when interpreted as the degrees
of the Reidemeister torsion of a certain twisted chain complex.
%The bounds are very powerful and
%can be easily implemented with a computer.
We show that these lower bounds give the correct genus bounds for all knots with 12 crossings or
less, including the Conway knot and the Kinoshita--Terasaka knot which have trivial Alexander
polynomial.
%\footnote{This will be the first of many occasions where I slightly change your wording,
%admittedly it sometimes comes down to taste. I guess we'll have to talk about this on the phone at
%some point. Your version was:\\
%`The genera of all knots with 12 crossings or less are known, and our computations show that the
%above new lower bounds give the sharp bounds for the genera of these knots.'}
%We also give many examples of closed manifolds and link
%complements where twisted Alexander polynomials detect the correct
%Thurston norm.

We also give obstructions to fibering 3--manifolds using twisted Alexander polynomials and detect
all knots with 12 crossings or less that are not fibered. For some of these it was unknown whether
or not they are fibered. Our work in particular extends the fibering obstructions of Cha to the
case of closed manifolds.
\end{abstract}
\maketitle

\section{Introduction}
%-------------------
\subsection{Definitions and history}
  Let $M$ be a 3--manifold. Throughout the paper we will assume that all
3--manifolds are compact, orientable and connected.  Let $\phi\in H^1(M)$ (integral coefficients
are understood). The \emph{Thurston norm} of $\phi$ is defined as
 \[
||\phi||_{T}=\min \{ \chi_-(S)\, | \, S \subset M \mbox{ properly embedded surface dual to }\phi\}.
\]
Here, given a surface $S$ with connected components $S_1\cup\dots \cup S_k$, we define
$\chi_-(S)=\sum_{i=1}^k \max\{-\chi(S_i),0\}$.
 Thurston
\cite{Th86} showed that this defines a seminorm on $H^1(M)$ which can be extended to a seminorm
on $H^1(M;\R)$. As an example consider $X(K):=S^3\sm \nu K$, where $K\subset S^3$ is a knot and
$\nu K$ denotes an open tubular  neighborhood of $K$ in $S^3$. Let $\phi\in H^1(X(K))$ be a
generator, then it is easy to see that $||\phi||_T=2\, \gen(K)-1$.

%It is an important problem to find methods for computing  the Thurston norm. Such methods  have
%many applications even outside of topology. For example using work of Freedman and He \cite{FH91} bounds
%on the Thurston norm translate into lower bounds for the ropelength
% \cite{CKS02} (cf. Section
%\ref{sectionrope}). Furthermore bounds on the Thurston norm also
%have applications in electrodynamics \cite{CK02,Ko04}.

It is a classical result of Alexander that
 \[ 2\mbox{genus}(K) \geq \deg(\Delta_K(t)),   \]
 where $\Delta_K(t)$ denotes the Alexander polynomial of a knot $K$. In recent years this was greatly
generalized. Let $M$ be a 3--manifold whose boundary is empty or
consists of tori. Let $\phi \in H^1(M) \cong \Hom(H_1(M),\Z)$ be
primitive, i.e., the corresponding homomorphism $\phi:H_1(M)\to
\Z$ is surjective. Then McMullen \cite{Mc02} showed that if the
Alexander polynomial $\alex\in \Q[t^{\pm 1}]$ of $(M,\phi)$ is
non--zero, then
 \[ ||\phi||_T \geq \deg\left(\alex\right)-(1+b_3(M)).\]
 Here $b_3(M)$ denotes the third Betti number of $M$, in particular  $b_3(M)=1$ if $M$ is closed and
$b_3(M)=0$ if $M$ has boundary.
% where $\ti{b}_3(M):=1$ if $M$ is closed
%or if all of its boundary components are spheres, and
%$\ti{b}_3(M):=0$ otherwise.
An alternative proof for closed manifolds was given by Vidussi \cite{Vi99,Vi03} using results of
Kronheimer--Mrowka \cite{KM97} and Meng--Taubes \cite{MT96} in Seiberg--Witten theory.
% We refer to
%\cite{Kr98,Kr99} for more on the connection between the Thurston norm, Seiberg--Witten theory and
%4--dimensional geometry.

Harvey \cite{Ha05} in the general case and Cochran \cite{Co04} in the knot complement case
generalized McMullen's inequality. They showed that the degrees of \emph{higher-order} Alexander
polynomials which are defined over non--commutative polynomial rings give lower bounds on the
Thurston norm. Later Harvey's work \cite{Ha05} was refined by Turaev \cite{Tu02b}.

%Cochran \cite{Co04} in the knot complement case and Harvey
%\cite{Ha05} and Turaev \cite{Tu02a,Tu02b} in the general case
%generalized McMullen's inequality. They studied maps
%$\Z[\pi_1(M)]\to \kt$ where $\K$ is a skew field and $\kt$ is a
%skew Laurent polynomial ring. They showed that the degrees of
%corresponding \emph{higher-order} Alexander polynomials give lower
%bounds on the Thurston norm.

In this paper we will show how the degrees of \emph{twisted Alexander polynomials} give lower
bounds on the Thurston norm.
%These bounds are easy to compute and remarkably strong.

%---------------------------------------
\subsection{Twisted Alexander polynomials and Reidemeister torsion}

In the following let $\F$ be  field. Let $\phi\in H^1(M) \cong \Hom(\pi_1(M),\Z)$ and
$\a:\pi_1(M)\to \gl(\F,k)$ a representation. Then $\a\otimes \phi$ induces an action of
$\pi_1(M)$ on $\F^k\otimes_\F \ft=:\fkt$ and we can therefore consider the twisted homology
$\ft$--module $H_i^{\a}(M;\fkt)$. We define $\Delta_{i}^\a(t) \in \ft$ to be its order; it is
called \emph{the $i$-th twisted Alexander polynomial of $(M,\phi, \a)$} and well--defined up to
multiplication by a unit in $\ft$. The twisted Alexander polynomial of a knot was introduced by
Lin \cite{Lin01} in 1990. In this paper we use the above homological definition of Kirk and
Livingston \cite{KL99}. These polynomials can be computed efficiently using Fox calculus and
Poincar\'e duality for twisted homology.  We refer to Section \ref{sectiondefalex} for more
details.

%In this paper we follow the definition of twisted Alexander
%polynomials given by Kirk and Livingston \cite{KL99}.

If $\partial M$ is empty or consists of tori and if $\twialex \ne 0$, then we will show that
$H_i^\a(M;\fkt \otimes_{\F\tpm} \F(t))=0$ for all $i$. Therefore the Reidemeister torsion
$\tau(M,\phi,\a)\in \F(t)$ is defined (cf. \cite{Tu01} for a definition) and (cf.
\cite[p.~20]{Tu01} or \cite[Theorem~3.4]{KL99})
\[ \tau(M,\phi,\a)=\prod_{i=0}^{2} \Delta_i^\a(t)^{(-1)^{i+1}}\in \F(t). \]
The equality holds up to multiplication by a unit in $\ft$.
% and will not make use of the fact that $\tau(M,\phi,\a)$ has
%in general a smaller indeterminacy.
For $f(t)/g(t)\in \f(t)$ we define $\deg(f(t)/g(t)):=\deg(f(t))-\deg(g(t))$ for $f(t),g(t)\in
\ft$. This allows us to consider $\deg(\tau(M,\phi,\a))$.

%---------------------------------------
\subsection{Lower bounds on the Thurston norm}

%Let $\phi \in H^1(M)$ be primitive and $\a:\pi_1(M)\to \glfk$ a
%representation. We write
%\[ \chi^\a_-(M):=-\sum_{i=0}^3(-1)^i\dim_{\F}(H_i^\a(M;\fkt))=-\sum_{i=0}^3(-1)^i\deg\left(\twialexi\right)\]
%if all homology groups are finite--dimensional over $\F$.
The following theorem shows that the degrees of twisted Alexander polynomials can be used to give
lower bounds on the Thurston norm.

%We equip $\F^k$ with the canonical hermitian inner product. We denote the group of unitary matrices
%of size $k$ by $\ufk$.

\begin{theorem} \label{mainthm} Let $M$ be a 3--manifold whose boundary is empty or consists of tori.
Let $\phi \in H^1(M)$ be non--trivial and $\a:\pi_1(M)\to \glfk$ a representation such that
$\twialex \ne 0$. Then
\[  \tnphi \geq  \frac{1}{k} \deg(\tau(M,\phi,\a)).\]
Equivalently,
\[ \tnphi \geq \frac{1}{k}\big(\degtwialex- \degtwialexzero -\degtwialextwo \big).\]
\end{theorem}

The proof of Theorem \ref{mainthm} is partly based on ideas of McMullen \cite{Mc02} and Turaev
\cite{Tu02b}. For one--dimensional representations it is easy to determine $\twialexzero$ and
$\twialextwo $ and one can easily show that Theorem \ref{mainthm}  contains McMullen's bound for
one--variable Alexander polynomials (\cite[Proposition 6.1]{Mc02}) and results of Turaev
\cite{Tu02a}.

In Theorem \ref{thmthurstonlink} we give lower bounds in the case $\twialex= 0$ for certain
$\phi$'s. In \cite{FK05} we introduce twisted Alexander norms (similar to McMullen's Alexander
norm \cite{Mc02}) which are well--suited to study the Thurston norm of link complements. We also
refer to \cite{Fr06} for a further extension of Theorem \ref{mainthm}.

%In \cite{Fr06} we will prove a version of Theorem \ref{mainthm}
%over skew fields, which combines our lower bounds from Theorem
%\ref{mainthm} with the lower bounds of Cochran, Harvey and Turaev
%\cite{Co04, Ha05, Tu02b}.

%We concentrate on proving Theorem
%\ref{mainthm}, i.e., the case for the modules over a commutative ring, and we only point out the
%changes to the proof of Theorem \ref{mainthm} which have to be made to prove the
%non--commutative generalization.

%An important source of representations is given by homomorphisms
%$\a:\pi_1(M)\to G$, $G$ a finite group. This induces a
%representation $\a:\pi_1(M)\to G\to \gl(\F,|G|)$ where the map $G
%\to \gl(\F, |G|)$ is the regular representation of $G$. (Note that
%$\gl(\F, |G|)$ is isomorphic to $\gl(\F[G])$.) In Section
%\ref{sectionfinitecovers} we give an elegant short proof of
%Theorem \ref{mainthm} in the case of a representation $\pi_1(M)\to
%G\to \gl(\F, |G|)$, using only McMullen's theorem and well--known
%properties of finite covers.

%---------------------------------------
\subsection{Fibered manifolds} \label{section:introfib} Let $\phi\in H^1(M)$ be non--trivial. We say
\emph{$(M,\phi)$ fibers over $S^{1}$} if the homotopy class of maps $M\to S^1$ induced by
$\phi:\pi_1(M)\to H_1(M)\to\Z$ contains a representative that is a fiber bundle over $S^{1}$. If
$K$ is a fibered knot, i.e., if $X(K)$ fibers, then it is a classical result of Neuwirth that
$K$ satisfies \be \label{abelianfib} \alexk \mbox{ is monic and } \deg(\alexk)=2 \, \gen(K).\ee

\begin{theorem} \label{mainthm2} Assume that $(M,\phi)$ fibers over
$S^1$ and that $M\ne S^1\times D^2, M\ne S^1\times S^2$. Let $\a:\pi_1(M)\to \glfk$ be a
representation. Then $\twialex\ne 0$ and
 \[ \tnphi = \frac{1}{k} \deg(\tau(M,\phi,\a)). \]
\end{theorem}

This result clearly generalizes the first classical condition on fibered knots. McMullen, Cochran,
Harvey and Turaev prove corresponding theorems in their respective papers \cite{Mc02, Co04, Ha05,
Tu02b}.
%Below we will show that Theorem \ref{mainthm2} gives a surprisingly strong fiberedness
%obstruction theorem.

Now let $R$ be a
Noetherian unique factorization domain (henceforth UFD), for example $R=\Z$. Given a representation
$\pi_1(M)\to \gl(R,k)$ Cha \cite{Ch03} defined a twisted Alexander polynomial $\twialex\in R\tpm$,
which is well--defined up to multiplication by a unit in $R\tpm$.
%This is a
%generalization of the Alexander polynomial $\Delta_K(t)\in \zt$ and coincides with the first
%twisted Alexander polynomial defined in Section \ref{sectiontwistedpoly} in the case that $R$ is a
%field.
%We say a polynomial $\twialex\in R\tpm$ is \emph{monic}, if its top coefficient is a unit in
%$R$.
Cha showed that for a fibered knot the polynomials $\Delta_1^\a(t)$ are monic \cite{Ch03}.
(Recall that a polynomial is called  monic, if its highest and lowest coefficient are units in
$R$.) Using Theorem \ref{mainthm2} we obtain the following theorem.
% that in particular generalizes Cha's obstructions to closed 3--manifolds.

\begin{theorem} \label{corfibered2}\label{mainthm3}
Let $M$ be a 3--manifold. Let $\phi\in H^1(M)$ be non--trivial such that  $(M,\phi)$ fibers over
$S^1$ and such that $M\ne S^1\times D^2, M\ne S^1\times S^2$. Let $R$ be a Noetherian UFD and let
$\a:\pi_1(M)\to \gl(R,k)$ be a representation.
%(In particular this is the case if $\a:\pi_1(M)\to G\to \gl(R[G])$ is a
%representation with $G$ a finite group.)
Then $\twialex \in R\tpm$ is monic and
\[ \tnphi =
\frac{1}{k} \deg(\tau(M,\phi,\a)). \]
\end{theorem}

In fact in Proposition \ref{prop:rewrite} we show that if the fibering obstruction of Theorem
\ref{mainthm2} vanishes, then the conclusion of Theorem \ref{corfibered2} holds.
%In particular
%this shows that our fibering obstruction of Theorem \ref{mainthm2} also contains the second
This shows that the obstructions of Theorem \ref{mainthm2} contain Neuwirth's and Cha's \cite{Ch03}
obstructions for fibered knots and extend them to closed 3--manifolds. Theorem \ref{mainthm2} is
also closely related to the result of \cite{GKM05} on fibered knots.
%Note that Theorem \ref{corfibered2} generalizes Cha's obstructions
%to closed 3--manifolds.

%---------------------------------------
\subsection{Examples}\label{introexample}
We give two main examples in Section \ref{sectionexamples}. First we show that the lower bounds of
Theorem \ref{mainthm} give, for appropriate representations, the correct genus bounds for all knots
with up to 12 crossings. These genera have  been found  by Gabai, Rasmussen, Stoimenow et. al. (cf.
\cite{CL} and \cite{Sto}).
%\footnote{Again, a matter of taste. I like my version since it's shorter
%and I don't think we have to mention every time that that was known. I think that what I wrote does
%not make any unfair claims.\\ "First, using Theorem ref{mainthmunitary} [we don't need the special
%unitary the special unitary theorem, since $\partial M\ne \emptyset$!!!!!], which is a special case
%of Theorem \ref{mainthm}, we compute lower bounds of the genera of knots with up to 12 crossings.
%The genera of these knots are already known (cf. \cite{CL}) and our computations show that the
%lower bounds from Theorem ref{mainthmunitary} give the right genus bounds for all these knots."}
Note that some knots with up to 12 crossings have trivial Alexander polynomial, and hence the genus
bounds of McMullen, Cochran and Harvey vanish. These, and all later computations in this paper,
were done using the program \emph{Knotwister} \cite{Fr05}.
%\footnote{I deleted this sentence:\\
%"For example, the Conway knot $K=11_{401}$ (\emph{knotscape} notation, cf. \cite{HT}) has
%$\Delta_K(t) = 1$ and the genus of the Conway knot is 3."\\
%If we want such a sentence, then we have to talk more how $\Delta^\a$ determines the genus,
%otherwise this sentence is pointless.}

Second we apply Theorem \ref{mainthm2}
 to study the fiberedness of knots. It is known that
a knot $K$ with 11 or fewer crossings is fibered if and only if $K$ satisfies Neuwirth's
condition (\ref{abelianfib}). Hirasawa and Stoimenow \cite{Sto} had started a program to find
all fibered 12--crossing knots. Using methods of Gabai they showed that except for thirteen
knots a 12--crossing knot is fibered if and only if it satisfies condition (\ref{abelianfib}).
Furthermore they showed that among these 13 knots the knots $12_{1498}$, $12_{1502}$,
$12_{1546}$ and $12_{1752}$ are not fibered even though they satisfy condition
(\ref{abelianfib}). Using Theorem \ref{mainthm2} we showed the non--fiberedness of these 4 knots
and we also showed that the remaining  9 knots are not fibered either. These 9 knots are:
$$
12_{1345},12_{1567},12_{1670},12_{1682},12_{1771},
12_{1823},12_{1938},12_{2089},12_{2103}.
$$
%To our knowledge these 9 knots are new. Stoimenow and Hirasawa showed that all the remaining
%12--crossing knots are fibered if and only if the Alexander polynomial is monic and $2\,
%\gen(K)=\deg(\alexk)$.
This result completes the classification of all fibered 12--crossing knots. We note that later
Jacob Rasmussen also showed that these 13 knots are not fibered using knot Floer homology (cf.
\cite[Section 3]{OS05}).

%As we pointed out our fibering obstructions work for closed
%manifolds as well. If $K$ is one of the 13 12--crossing knots in
%the previous paragraph, then we can easily show using Theorem
%\ref{mainthm2} and \emph{KnotTwister} that the zero surgery on $K$
%in $S^3$ is not fibered. (See Section \ref{sectionexamples11}.)
%
%The situation for links is more complex. On the one hand  in many
%interesting cases twisted Alexander polynomials give the correct
%bound. For example in Section \ref{sectionrope} we reprove results
%of Harvey on the ropelength of a certain link \cite{Ha05}. We also
%successfully apply our theory in Section \ref{section:sat}. On the
%other hand boundary links have mostly vanishing twisted Alexander
%polynomials and therefore our lower bounds do not apply in
%general. But in Section \ref{sectiontorsion} we show that in some
%cases we can still extract lower bounds from the degrees of
%twisted Alexander polynomials corresponding to the $\ft$--torsion
%submodule of $H_1^\a(X(L);\fkt)$ where $X(L)$ is the link
%complement in the 3--sphere (cf. Theorem \ref{thmthurstonlink}).

%---------------------------------------
\subsection{Outline of the paper}
In Section \ref{sectionthurstonalex} we give a definition of twisted Alexander polynomials and
we discuss the Alexander polynomials of 3--manifolds.  We give the proofs of Theorem
\ref{mainthm} and Theorem \ref{mainthm2} in Section \ref{sectionproofmainthm}. In Section
\ref{sectiontorsion} we discuss the case that  $\twialex=0$. We discuss the examples in Section
\ref{sectionexamples}. In Section \ref{sectioncha} we prove Theorem \ref{mainthm3}.
% Finally in
%Section \ref{sectionconj} we discuss and give further evidence for
%Conjectures \ref{conjcha} and related conjectures.
\\

\noindent {\bf Notations and conventions:} We assume that all 3--manifolds are compact, oriented
and connected. All homology groups and all cohomology groups are with respect to
$\Z$--coefficients, unless it specifically says otherwise. For a link $L$ in $S^3$, $X(L)$
denotes the exterior of $L$ in $S^3$. (That is, $X(L) = S^3\sm \nu L$ where $\nu L$ is an open
tubular neighborhood of $L$ in $S^3$). $\F$ will always denote a  field. We identify the group
ring $\F[\Z]$ with $\ft$. For a 3--manifold $M$ we use the canonical isomorphisms to identify
$H^1(M) = \Hom(H_1(M), \Z) = \Hom(\pi_1(M),\Z)$. Hence sometimes $\phi\in H^1(M)$ is regarded as
a homomorphism $\phi : \pi_1(M) \to \Z$ (or $\phi : H_1(M) \to \Z$) depending on the context.
\\

\noindent {\bf Acknowledgments:} The authors would like to thank Alexander Stoimenow for providing
braid descriptions for the examples and Stefano Vidussi for pointing out the advantages of using
Reidemeister torsion. The first author would also like to thank Jerry Levine for helpful
discussions and he is indebted to Alexander Stoimenow for important feedback on the program
\emph{KnotTwister}.
%This work was supported by the faculty research fund of Konkuk University in
%2006

%Here
%$\gen(K)$ is defined to be the minimal genus of a surface bounding
%$K$. Using Haken's theory of normal surfaces Schubert \cite{Sc61}
%gave decision procedures on finding the genus of a knot, but these are not practical to implement.

%===========================================
\section{The  twisted Alexander polynomials and duality}\label{sectionthurstonalex} \label{sectiondefalex}

%==============================================
\subsection{Twisted homology groups} \label{section:twihom}
We first give a   definition of twisted homology groups and discuss some of their properties.
Let $X$ be a topological space, $Y\subset X$ a (possibly empty) subset and $x_0\in X$ a point.
Let $R$ be a ring (e.g. $R=\F$ or $R=\ft$) and $\b:\pi_1(X,x_0)\to \gl(R,k)$ a representation.
This naturally induces a left $\Z[\pi_1(X,x_0)]$--module structure on $R^k$.

Denote by $\ti{X}$ the set of all homotopy classes of paths starting at $x_0$ with the usual
topology. Then the evaluation map $p:\ti{X}\to X$ is the universal cover of $X$. Note that $g\in
\pi_1(X,x_0)$ naturally acts on $\ti{X}$ on the right by precomposing any path by $g^{-1}$.

Given $Y\subset X$ we let $\ti{Y}=p^{-1}(Y)\subset \ti{X}$. Then the above $\pi_1(X,x_0)$ action on
$\ti{X}$ gives rise  to a right $\Z[\pi_1(X,x_0)]$--module structure on the chain groups
$C_*(\ti{X}),C_*(\ti{Y})$ and $C_*(\ti{X},\ti{Y})$. Therefore we can form the tensor product over
$\Z[\pi_1(X,x_0)]$ with $R^k$, we define
\[ \ba{rcl}
H_i^{\b}(X;R^k)&=&H_i(C_*(\ti{X})\otimes_{\Z[\pi_1(X,x_0)]}R^k)),\\
H_i^\b(Y\subset X;R^k)&=&H_i(C_*(\ti{Y})\otimes_{\Z[\pi_1(X,x_0)]}R^k),\\
H_i^{\b}(X,Y;R^k)&=&H_i(C_*(\ti{X},\ti{Y})\otimes_{\Z[\pi_1(X,x_0)]}R^k)). \ea\] Sometimes we
refer to $H_i^\b(Y\subset X;R^k)$ as \emph{twisted subspace homology}. Note that if we have
inclusions $Z\subset Y\subset X$ then we get an induced map $H_i^\b(Z\subset X;R^k)\to
H_i^\b(Y\subset X;R^k)$. Also note that we have an exact sequence of complexes
\[ 0 \to C_i(\ti{Y})\otimes_{\Z[\pi_1(X,x_0)]}R^k\to
C_i(\ti{X})\otimes_{\Z[\pi_1(X,x_0)]}R^k\to C_i(\ti{X},\ti{Y})\otimes_{\Z[\pi_1(X,x_0)]}R^k\to
0\] which gives rise to a long exact sequence \be \label{longexact} \dots \to H_i^{\b}(Y\subset
X;R^k)\to H_i^{\b}(X;R^k)\to H_i^{\b}(X,Y;R^k)\to \dots .\ee

Now denote by $Y_i,i\in I,$ the path connected components of $Y$. Pick base points $y_i\in Y_i,i\in
I,$ and paths $\g_i:[0,1]\to Y$ with $\g_i(0)=y_i$ and $\g_i(1)=x_0$. Then we can get induced
representations $\b_i(\g_i):\pi_1(Y_i,y_i)\to \pi_1(X,y_i)\to \pi_1(X,x_0)\to \gl(R,k)$ and induced
homology groups $H_j^{\b_i(\g_i)}(Y_i;R^k)$ using the universal cover of $Y_i$.

\begin{lemma} \label{lem:twiiso}
Given $\g_i$ there exists a canonical  isomorphism
\[ H_j^{\b_i}(Y_i\subset X;R^k)\cong  H_j^{\b_i(\g_i)}(Y_i;R^k).\]
\end{lemma}

\begin{proof}
Let $K$ be the image of $\pi_1(Y_i,y_i)$ under the map $\i(\g_i):\pi_1(Y_i,y_i)\to \pi_1(X,y_i)\to
\pi_1(X,x_0)$ induced by $\g_i$. Denote by $\widetilde{Y_i}^{K}$ the  cover of $Y_i$ corresponding
to $\pi_1(Y_i,y_i)\to K$. More precisely, we take
\[  \widetilde{Y_i}^{K}=\{ \s:[0,1]\to Y_i | \s(0)=y_i \}/ \sim \]
where $\sim$ is the equivalence relation given by
\[  \s_1 \sim \s_2  \mbox{ if } \s_1(1)=\s_2(1) \mbox{ and }\i(\g_i)(\s_1\s_2^{-1})=e\in K.\]
 Then we get a
well--defined injective map
%\footnote{I don't think we have to, or even should define $\ti{Y_i}$,
%we already defined it for any subset of $X$. I adopted the convention that for paths $a,b$ the path
%$ab$ means, $a$ first, $b$ second, which I think is the more common convention.}
\[ \ba{rcl} \widetilde{Y_i}^{K}&\to& \ti{Y_i}\subset \ti{X}\\
  \,   [\s]&\mapsto& [\g_i^{-1}\s].\ea \] We will use this injection to identify
 $\widetilde{Y_i}^{K}$ with its image in $\widetilde{Y_i}$. Note that $\widetilde{Y_i}$ is the
disjoint union of copies of $\widetilde{Y_i}^K$ indexed by $\pi_1(X,x_0)/K$.
 In particular a singular simplex in $\widetilde{Y_i}$ is of
the form $\s g$ for a singular simplex $\s$ in $\widetilde{Y_i}^{K}$ and an element $g\in
\pi_1(X,x_0)$. Mapping
\[ \s g\otimes_{\Z[\pi_1(X,x_0)]}  v \mapsto \s\otimes_{\Z[K]} \b(g)v\]
(with $v\in R^k)$ induces an  isomorphism $C_*(\widetilde{Y_i})\otimes_{\Z[\pi_1(X,x_0)]}R^k\cong
C_*(\widetilde{Y_i}^{K})\otimes_{\Z[K]}R^k$ of chain complexes. Denote the universal cover of $Y_i$
by $\widetilde{Y_i}^{\pi_1(Y_i,y_i)}$. It is easy to see that the projection induced map
$C_*(\widetilde{Y_i}^{\pi_1(Y_i,y_i)})\to C_*(\widetilde{Y_i}^{K})$ gives rise to  an isomorphism
of chain complexes:
\[
C_*(\widetilde{Y_i}^{K})\otimes_{\Z[K]}R^k\cong
C_*(\widetilde{Y_i}^{\pi_1(Y_i,y_i)})\otimes_{\Z[\pi_1(Y_i,y_i)]}R^k.\]
\end{proof}

Note that the isomorphism of the lemma only depends on the choice of $\g_i$, and we call the
isomorphism the \emph{isomorphism induced by $\g_i$}.

It is clear that the isomorphism type of $H_j^\b(X,Y;R^k)$ does not depend on the choice of the
base point. In most situations we can and will therefore suppress the base point in the notation
and the arguments. We will also normally write $\b$ instead of $\b(\g_i)$. Furthermore we write
$H_j^\b(Y;R^k)=\oplus_{i\in I}H_j^\b(Y_i;R^k)$. With these conventions  the long exact sequence
(\ref{longexact}) induces a long exact sequence
\[ \dots \to H_j^\b(Y;R^k)\to H_j^\b(X;R^k) \to H_j^\b(X,Y;R^k)\to \dots.\]
The isomorphism type of the sequence  is independent of all the choices made.

%For the remainder of this section we can safely stop mentioning the base points and the connecting
%paths since the isomorphism class of the objects and the maps are independent. But when we
%formulate and proof Proposition \ref{prop:exact} we will need to argue carefully with base points
%and connecting paths.

%==============================================
\subsection{The twisted Alexander
polynomials}\label{sectiontwistedpoly}

For the remainder of this section we assume that  $M$ is a 3--manifold and $\phi\in H^1(M)$. Let
$\a:\pi_1(M)\to \gl(\F,k)$ be a representation. We can now define a left $\Z[\pi_1(M)]$--module
structure on $\F^k\otimes_\F \ft=:\fkt$ via $\a\otimes \phi$ as follows:
\[  g\cdot (v\otimes p):= (\a(g)\cdot v)\otimes (\phi(g)\cdot p) = (\a(g) \cdot v)\otimes (t^{\phi(g)}p) \]
where $g\in \pi_1(M), v\otimes p \in \F^k\otimes_\F \ft = \fkt$. Put differently, we get a
representation $\a\otimes \phi:\pi_1(M)\to \gl(\ft,k)$. We call $H_i^{\a\otimes \phi}(M;\fkt)$ the
 \emph{$i$--th twisted Alexander module of $(M,\phi,\a)$}.
Usually we drop the notation $\phi$ and write $H_*^{\a}(M;\fkt)$.
Note that $H_i^{\a}(M;\fkt)$ is a finitely generated module over
the PID $\ft$. Therefore there exists an isomorphism
 \[ H_i^{\a}(M;\fkt) \cong  \ft^f \oplus \bigoplus_{i=1}^l \ft/(p_i(t))
\] for $p_1(t),\dots,p_l(t) \in \ft\sm \{0\}$.
 We define
\[ \Delta^\a_{M,\phi,i}(t):=\left\{ \ba{rl} \prod_{i=1}^l p_i(t), &\mbox{ if }f=0\\
0,&\mbox{ if }f>0.\ea\right. \]
 This is called the \emph{$i$--th twisted Alexander polynomial} of
$(M,\phi,\a)$. We furthermore define $\ti{\Delta}^\a_{M,\phi,i}(t):=\prod_{i=1}^k p_i(t)$
regardless of $f$. In most cases we drop the notations $M$ and $\phi$ and write $\twialexi$ and
$\titwialexi$. It follows from the structure theorem of finitely generated modules over a PID
that these polynomials are well--defined up to multiplication by a unit in $\ft$.
%In Section \ref{sectionfox} we will see that $\twialexi$
%and $\titwialexi$ can be computed easily for $i=0,1$ given a
%presentation of $\pi_1(M)$.

%\begin{remark}
%%The twisted Alexander polynomial for a knot was originally defined by Lin in 1990 using a
%%presentation of the fundamental group \cite{Lin01}. This was generalized by Jiang and Wang
%%\cite{JW93} and  by Wada \cite{Wa94} given only a presentation of a group and a representation
%%to $\gl(R,k)$ where $R$ is a UFD. Note that Wada's definition differs slightly from our
%%definition. Our homological definition of twisted Alexander polynomials in the above was
%%originally introduced by Kirk and Livingston in \cite{KL99}.
%The twisted Alexander polynomial of a knot was introduced by Lin
%\cite{Lin01} in 1990. Various versions of twisted Alexander
%polynomials have been successfully used in many situations to
%provide more information than can be extracted from the untwisted
%Alexander polynomial
%\cite{JW93,Wa94,Kit96,KL99,KL99b,Ch03,HLN04}. In particular we
%note that Kirk and Livingston \cite{KL99} first introduced the
%above homological definition of twisted Alexander polynomials for
%a finite complex. We refer to \cite[Section 4]{KL99} for the
%relationship between our definition and the other definitions of
%twisted Alexander polynomials.
%\end{remark}

For an oriented knot $K$ we always assume that  $\phi$ denotes the generator of $H^1(X(K))$
given by the orientation. If $\a: \pi_1(X(K))\to \gl(\Q,1)$ is the trivial representation then
the Alexander polynomial $\Delta_1^{\a}(t)$ equals  the classical Alexander polynomial
$\Delta_K(t) \in \Q[t^{\pm 1}]$ of the knot $K$.

%\begin{remark}
%When $K$ is a knot, then the untwisted homology $H_1(X(K);\ft)$ is
%$\ft$--torsion. But even for a knot complement it can happen that
%in the twisted case $H_1^\a(X(K);\fkt)$ is not $\ft$--torsion (cf.
%e.g. \cite{KL99}).
%\end{remark}

Let $f=\sum_{i=m}^n a_it^i\in \ft\sm \{0\}$ with $a_m\ne 0, a_n \ne 0$. Then we define
$\deg(f)=n-m$. The following observation follows immediately from the classification theorem of
finitely generated modules over a PID.

\begin{lemma} \label{lemmafinite}
$\twihomi$ is a finite--dimensional $\F$--vector space if and only
if $\twialexi\ne 0$. If $\twialexi \ne 0$, then
\[ \degtwialexi=\dimm_{\F}\left(\twihomi\right). \]
 Furthermore $\deg(\titwialexi)=\dimm_{\F}\left(\mbox{Tor}_{\ft}(\twihomi)\right)$.
\end{lemma}

%%===========================================
\subsection{Duality for twisted homology} \label{sectionduality}
In this section we discuss a
 duality theorem for twisted homology which we will need to compute higher
twisted Alexander polynomials of 3--manifolds and which  will also play an important role in the
proof of Proposition  \ref{prop:b2}.

Let $\F$ be a field with (possibly trivial) involution $f\mapsto \ol{f}$.  We equip $\fk$ with the
standard hermitian inner product $\langle v ,w \rangle=v^t\ol{w}$ (where we view elements in $\fk$
as column vectors). We extend the involution on $\F$ to $\ft$ by taking $t \mapsto t^{-1}$. We
equip $\fkt$ with the hermitian inner product defined by $\langle vt^i,wt^j\rangle:=\langle
v,w\rangle t^it^{-j}$ for all $v,w\in \fk$.

In the following let $R=\F$ or $R=\ft$. Let $\b:\pi_1(M)\to \gl(R,k)$ be a representation. There
exists a unique representation $\ol{\b} : \pi_1(M) \to \gl(R,k)$  such that
\[ \langle \b(g^{-1})v,w\rangle =\langle v,\ol{\b}(g)w\rangle \]
for all $v,w\in \F^k,g\in \pi_1(M)$.

%The following duality theorem for homology with twisted coefficients is crucial for computing
%$\Delta^{\a}_2$ and in the proof of Proposition \ref{prop:b2}.

The following Lemma is a variation on \cite[p.~639]{KL99}.

\begin{lemma} \label{lemmadualitybig}
Let $X$ be an $n$--manifold and $\b:\pi_1(X)\to \gl(R,k)$ a representation. Then
\[
H_{n-i}^\b(X;R^k)\cong \hom_R(H_{i}^{\overline{\b}}(X,\partial X;R^k),R)\oplus
\ext_R(H_{i-1}^{\overline{\b}}(X,\partial X;R^k),R)
\] as
 $R$--modules.
\end{lemma}

%\noindent Lemma \ref{lemmadualitybig} easily follows from
%Poincar\'e duality and the universal coefficient theorem, and the
%proof is omitted.

\begin{proof}
%Let $V'=\overline{V}$ as $R$--modules equipped with the right
%$\Z[\pi_1(M)]$--module given by $v\cdot g:=\ol{\b}(g^{-1})v$ for
%$v\in \overline{V}$ and $g\in \pi_1(M)$.
Let $\pi:=\pi_1(X)$. We write $(R^k)'$ when we think of $R^k$ as equipped with the right
$\Z[\pi]$--module structure given by $v\cdot g:=\b(g^{-1})v$ for $v\in R^k$ and $g\in \pi$. By
Poincar\'e duality we have (recall that $\widetilde{\partial X}$ is the preimage of $\partial X$
under the covering map $\ti{X}\to X$)
\[ H_{n-i}^\b(X;R^k)\cong H^{i}_{\b}(X,\partial X;(R^k)')
:= H_i\big(\hom_{\Z[\pi]}(C_*(\ti{X},\widetilde{\partial X}),(R^k)')\big).\]  Using the inner
product we get an isomorphism of $R$--module chain complexes:
\[
\ba{rcl} \hom_{\Z[\pi]}(C_*(\ti{X},\widetilde{\partial
X}),(R^k)')&\to&
\hom_R\big(C_*^{\ol{\b}}(\ti{X},\widetilde{\partial X};
R^k),R\big)
=\hom_R\big(C_*(\ti{X},\widetilde{\partial X})\otimes_{\Z[\pi]} R^k,R\big)\\
f&\mapsto& \left((c\otimes w)\mapsto \langle f(c),w\rangle\right). \ea
\]
%Note that this map is well--defined since $\langle \b(g^{-1})v,w\rangle =\langle
%v,\ol{\b}(g)w\rangle$. It is now easy to see that this defines in fact an isomorphism of
%$R$--module chain complexes.
The lemma now follows from applying the universal coefficient theorem.
% for chain complexes over the
%PID $R$ to $C_*(\ti{X},\partial \ti{X})\otimes_{\Z[\pi]}V(\ol{\b})$. The lemma is now immediate.
\end{proof}

%===========================================
\subsection{Twisted Alexander polynomials of 3--manifolds}

\begin{lemma} \label{lem:delta03}
Let $\phi \in H^1(M)$ be non--trivial and $\a:\pi_1(M)\to \glfk$ a representation. Then
\bn
\item $\Delta_0^\a(t)\ne 0$,
\item $\Delta_3^\a(t)=1$.
\en
\end{lemma}

\begin{proof}
Both statements follow from a straightforward  argument using a cell decomposition of $M$ as in the
proof of Proposition \ref{prop:rewrite}.
%\footnote{There's also another alternative:
%$H_i(M;\fkt)=H_i(\hat{M};\fk)$, where $\hat{M}$ is the infinite cyclic cover. Then the lemma
%follows easily from the fact that $\hat{M}$ is homotopy equivalent to a 2--complex and that
%$\hat{M}$ has only finitely many components. Which proof do you think should we use?}
Alternatively note that Kirk and Livingston showed (1) in \cite[Proposition 3.5]{KL99}. For (2) we
apply Lemma \ref{lemmadualitybig} with $R= \ft$ and $\b = \a\otimes \phi$, and get
 \[ H_3^{\a\otimes \phi}(M;\fkt)\cong \hom_{\ft}\left(H_0^{\overline{\a\otimes \phi}}(M,\partial
 M;\fkt),\ft\right)
\] as $\ft$--modules. Note that $\overline{\a\otimes \phi}=\ol{\a}\otimes (-\phi)$.
It follows from (1) that $H_0^{\ol{\a}\otimes (-\phi)}(M;\fkt)$ is
$\ft$--torsion. It follows from the long exact homology sequence
that $H_0^{\ol{\a}\otimes (-\phi)}(M,\partial
 M;\fkt)$ is $\ft$--torsion as well, hence $H_3^{\a\otimes \phi}(M;\fkt)=0$.
\end{proof}

\begin{proposition} \label{prop:delta02}  \label{prop:delta2}
Let $M$ be a 3--manifold whose boundary is empty or consists of tori and let $\phi \in H^1(M)$
be non--trivial. Let $\a:\pi_1(M)\to \glfk$ be a representation such that $\twialex \ne 0$.
 \bn
\item If $M$ is closed, then
%$\Delta_1^{\at}(t)\ne 0$ and
\[ \twialextwo=\Delta_0^{\at}(t^{-1}). \]
\item If $M$ has non--empty boundary, then $\twialextwo=1$. \en In
particular $\degtwialextwo =b_3(M)
\deg\big(\Delta_0^{\at}(t)\big)$. Furthermore, if $\a$ is unitary,
i.e. $\a=\at$, then $\degtwialextwo =b_3(M)\degtwialexzero$.
\end{proposition}

For the proof we need the following two useful lemmas which we will also need several times
later.

\begin{lemma} \label{lemmah0surjects}
Let $R$ be a ring, $A$  a group and $\a:A\to \gl(R,k)$ a representation. If $\varphi:B\to A$ is a
homomorphism, then $H_0^{\a \circ \varphi}(B;R^k)\to H_0^\a(A;R^k)$ is surjective. Furthermore if
$\varphi$ is an epimorphism, then $H_0^{\a \circ \varphi}(B;R^k)\to H_0^\a(A;R^k)$ is an
isomorphism.
\end{lemma}

\begin{proof}
The lemma follows immediately from the commutative diagram of exact sequences
\[ \ba{cccccccccc}
0&\to&\{ \a(\varphi(b))v-v | b\in B, v\in R^k\}&\to&R^k&\to& H_0^{\a \circ\v}(B;R^k)&\to&0\\
&&\downarrow&&\downarrow&&\downarrow \\
0&\to&\{ \a(a)v-v | a\in A, v\in R^k\}&\to&R^k&\to& H_0^\a(A;R^k)&\to&0 \ea
\] and the observation that the vertical map on the left is
injective (respectively an isomorphism).
\end{proof}

A standard argument shows the following  lemma.

\begin{lemma} \label{lemmab1} \label{lem:euler}
Let $X$ be an $n$--manifold, $\K$ a field (e.g. $\F$ or $\F(t)$),  and $\a:\pi_1(X)\to
\gl(\K,k)$ a representation. Then
\[ \sum_{i=0}^n (-1)^i \dim_{\K} (H_*^{\a}(X;\K^k))=k\chi(X). \]
\end{lemma}

\begin{proof}[Proof of Proposition  \ref{prop:delta02}]
We will first show that  $H_2^\a(M;\fkt)$ is $\ft$--torsion.  Note that it follows from the long
exact homology sequence for $(M,\partial M)$ and from duality that
$\chi(M)=\frac{1}{2}\chi(\partial M)$. Hence $\chi(M)=0$ in our case. It now follows from Lemma
\ref{lem:euler} (applied to the field $\f(t))$ that
\[ \sum_{i=0}^3 (-1)^i \dim_{\F(t)}\left( H_i^\a(M;\fkt\otimes_{\ft} \f(t))\right)=k\cdot\chi(M)=0.\]
Note that $H_i^\a(M;\fkt\otimes_{\ft} \f(t))=H_i^\a(M;\fkt)\otimes_{\ft} \f(t)$ since $\f(t)$ is
flat over $\ft$. By assumption and by Lemma \ref{lem:delta03} we have $H_i^\a(M;\fkt)\otimes_{\ft}
\f(t)=0$ for $i\ne 2$, hence $H_2^\a(M;\fkt)\otimes_{\ft} \f(t)=0$ as well.

Now we apply Lemma \ref{lemmadualitybig} and using that $\overline{\a\otimes \phi} = \at \otimes
(-\phi)$ we get
 \[ \ba{rccl} H_2^{\a\otimes \phi}(M;\fkt)&\cong & &\hom_{\ft}\left(H_1^{\at \otimes
(-\phi)}(M,\partial
 M;\fkt),\ft\right)\\
&&\oplus &\ext_{\ft}\big(H_0^{\at \otimes
(-\phi)}(M,\partial M;\fkt),\ft\big)\\
% &\cong &
%\overline{\hom_{\ft}\left(H_1^{\at}(M,\partial M;\fkt);\ft\right)} \oplus
%\overline{H_0^{\at}(M,\partial M;\fkt)}

\ea
\] as $\ft$--modules. Since we know that
$H_2^{\a\otimes \phi}(M;\fkt)$ is $\ft$--torsion it follows that the first summand on the right
hand side is zero.

By Lemma \ref{lem:delta03} $H_0^{\at \otimes (-\phi)}(M;\fkt)$ is $\ft$--torsion. From the long
exact homology sequence of the pair $(M,\partial M)$ it follows that $H_0^{\at \otimes
(-\phi)}(M,\partial M;\fkt)$ is also $\ft$--torsion. Since $H_0^{\at \otimes (-\phi)}(M,\partial
M;\fkt)$ is a finitely generated $\ft$--torsion module and $\ft$ is a PID,
$\ext_{\ft}(H_0^{\overline{\a\otimes \phi}}(M,\partial M;\fkt),\ft)\cong H_0^{\at \otimes
(-\phi)}(M,\partial M;\fkt)$.

If $M$ is closed then we get $ H_2^\a(M;\fkt)\cong H_0^{\at \otimes (-\phi)}(M;\fkt)$. Therefore we
deduce that $ \twialextwo=\Delta_0^{\at}(t^{-1}) $.
%Also note that
%$\hom_{\ft}\left(H_1^{\at}(M;\fkt);\ft\right)=0$ implies $\Delta_1^{\at}(t)\ne 0$.
If $\partial M\ne \emptyset$, then by Lemma \ref{lemmah0surjects} the map $H_0^{\at \otimes
(-\phi)}(\partial M;\fkt)\to H_0^{\at \otimes (-\phi)}(M;\fkt)$ is surjective, hence $H_0^{\at
\otimes (-\phi)}(M,\partial M;\fkt)=0$. This shows that $H_2^\a(M;\fkt)=0$ and hence
$\twialextwo=1$.
\end{proof}

\begin{remark}
Given a presentation for $\pi_1(M)$ the polynomials $\twialex$ and $\twialexzero$ can be computed
efficiently using Fox calculus (cf. e.g. \cite[p.~98]{CF77}, \cite{KL99}). We point out that
because we view $C_*(\ti{M})$ as a \emph{right} module over $\Z[\pi_1(M)]$ we need a slightly
different definition of Fox derivatives. We refer to \cite[Section~6]{Ha05} for details.
Proposition \ref{prop:delta02}  allows us to compute $\twialextwo$ using the algorithm for
computing the 0-th twisted Alexander polynomial.
\end{remark}

\begin{remark}
Let $\a:\pi_1(M)\to \glfk$ be a unitary representation. Then the inequality in Theorem
\ref{mainthm} becomes the computationally slightly simpler inequality
\[ \tnphi \geq
\frac{1}{k}\big(\degtwialex- \oneplusbm \degtwialexzero  \big).\]
\end{remark}

%===========================================
\subsection{Reidemeister torsion of 3--manifolds}

Assume that $\partial M$ is empty or consists of tori and that $\twialex \ne 0$. Then it follows
from Lemma \ref{lem:delta03} and Proposition \ref{prop:delta2} that $\twialexi \ne 0$ for all $i$
and hence $H_i^\a(M;\fkt \otimes_{\F\tpm} \F(t))=H_i^\a(M;\fkt )\otimes_{\F\tpm} \F(t)=0$ for all
$i$.  Therefore the Reidemeister torsion $\tau(M,\phi,\a)\in \F(t)^*/\{rt^l |r\in \F^*, l\in \Z\}$
is defined. We refer to \cite{Tu01} for an excellent introduction into the theory of Reidemeister
torsion.
% $\tau(M,\phi,\a)\in \F(t)$ is well--defined up to multiplication
%by an element of the form $rt^k, r\in \im\{\pi_1(M)\xrightarrow{\a} \glfk\xrightarrow{det}\F\}$.
%We will not make use of this, and

The following lemma follows from \cite[p.~20]{Tu01} combined with the fact that $\Delta^\a_3(t)
= 1$ (cf. also \cite[Theorem~3.4]{KL99})

\begin{lemma}  \label{lemmareidemeister}
If $\twialex\ne 0$, then $\tau(M,\phi,\a)$ is defined and
\[ \tau(M,\phi,\a)=\prod_{i=0}^{2} \Delta_i^\a(t)^{(-1)^{i+1}}\in \F(t) \]
 up to multiplication by a unit in
$\ft$.
\end{lemma}

For our purposes we can also  use this equality as a definition for $\tau(M,\phi,\a)$.
 We will mostly use $\tau(M,\phi,\a)$ as a convenient and concise way to
store information. We point out that $\tau(M,\phi,\a)$ can also be computed directly from the chain
complex of $M$ (cf. \cite{Tu01}).

%===========================================
\section{Proof of Theorems \ref{mainthm} and \ref{mainthm2}}\label{sectionproofmainthm}

%===========================================
\subsection{Proof of Theorem \ref{mainthm}}
For the remainder of this section let $M$ be a 3--manifold and let $\phi \in H^1(M)$ be primitive.
A \emph{weighted surface} $\hat{S}$ in $M$ is defined to be a collection of pairs $(S_i,w_i),
i=1,\dots,l$ where $S_i\subset M$ are properly disjointly embedded, oriented surfaces in $M$ and
$w_i$ are positive integers. We denote the union $\bigcup_{i} S_i \subset M$ by $S'$.

Every weighted surface $\hat{S}$ defines an element $\phi_{\hat{S}}:=\sum_{i=1}^l w_i\cdot
PD([S_i]) \in H^1(M)$ where $PD(f)\in H^1(M)$ denotes the Poincar\'e dual of an element $f\in
H_2(M,\partial M)$. By taking $w_i$ parallel copies of $S_i$ we get an (unweighted) properly
embedded oriented surface $S$ such that $\phi_{\hat{S}}=PD([S])$.
%An example of a the surface $\hat{S}^{\#}$ for a weighted surface
%$\hat{S}$ is given in
%Figure \ref{fig:weightedsurface}.
%
%\begin{figure}[h] \begin{center}
% \includegraphics[scale=0.35]{weightedsurface.eps} \caption{Weighted surface in a handlebody.}
%\label{fig:weightedsurface}
%\end{center}
% \end{figure}
We need the following very useful proposition proved by Turaev in \cite{Tu02b}.

\begin{proposition} \label{propminus}
 There exists a weighted surface $\hat{S}=(S_i,w_i)_{i=1,\dots,l}$ with
\bn
\item $\phi_{\hat{S}}=\phi$,
\item $\chi_-(S)=||\phi||_T$, and
\item $M\sm S'$ connected.
%\item $w_i\ne 0$ for all $i$.
%\item if $T\subset \partial M$ is a torus and if $\phi$
%vanishes on $H_1(T)$, then $T\cap \hat{S}^\#=\emptyset$.
\en
\end{proposition}

For the remainder of this section let $\hat{S}=(S_i,w_i)_{i=1,\dots,l}$ be a weighted surface as in
Proposition \ref{propminus}. We now do the following:
\bn
\item We pick orientation preserving disjoint
embeddings $\iota:S_i\times [0,w_i] \to M,i=1,\dots,l$ such that $\iota$ restricted to $ S_i\times
0$ is the identity (where we identify $S_i\times 0$ with $S_i$). We identify the image of $\iota$
with $ S_i \times[0,w_i]$.
\item Note that with these conventions we have $S=\cup_{i=1}^l \cup_{j=0}^{w_i-1} S_i\times j$ and $S'=\cup_{i=1}^l  S_i$.
\item For any subset $I\subset [0,1]$ we write $S\times I=\cup_{i=1}^l\cup_{j=0}^{w_i-1}S_i\times (j+I)$.
\item We let $\eps=\frac{1}{2}$.
\item We let $N=M\sm S\times (0,\eps)$
and we let $N'=M\sm \cup_{i=1}^l S_i\times (0,w_i-1+\eps)$. Note that $N'$ is connected by
Proposition \ref{propminus} (3).
%\item We denote the components of $N$ by $N_1,\dots,N_n$.
\item We pick a base point $m_0$ for $N'$ which also serves as a base point for $M$.
\item For $i=1,\dots,l$ pick a base point $s_i$ of $S_i$ and we pick a path in $N'$ connecting $m_0$ to $s_i$.
%\item To a subset $S_i\times [a,b]\subset M$ associate the base point $p_i\times a$ and connect it
%to $p$ using the obvious extension of the path from $p_i$ to $p$.
\en

Recall that given a representation $\b:\pi_1(M,m_0)\to \gl(R,k)$ we get, using the paths chosen
above, induced representations (and hence twisted homology groups) for $N'$ and $S_i\times
0,i=1,\dots,l$ which we denote by the same symbol. We will sometimes use the  isomorphisms of
Lemma \ref{lem:twiiso} induced by the chosen paths to identify the twisted homology groups with
the twisted subspace homology groups.

%Note that with these choices we get (using Lemma \ref{lem:twiiso}) inclusion induced maps
%\[ H_i^\b(S\times \eps;R^k)\cong H_i^\b(S\times \eps\subset N;R^k)\to H_i^\b(N;R^k)\]
% We also get a map
%\[ H_i^\b(N;R^k)\cong H_i^\b(N\subset M;R^k)\to H_i^\b(M;R^k).\]

In the following let $p:\ti{M}\to M$ be the universal cover of $M$ corresponding to the base
point $m_0$ as in Section \ref{section:twihom}, in particular $\ti{M}$ is the set of homotopy
classes of paths in $M$ starting at $m_0$. Also we again write $\ti{X}=p^{-1}(X)\subset \ti{M}$
for any $X\subset M$.
 Now note that given $a,a+\delta \in [0,w_i]$ we get
a $\Z[\pi_1(M,m_0)]$--equivariant map $f_{\delta}:p^{-1}(S_i\times a)\to p^{-1}(S_i\times
(a+\delta))$ by extending a path from $m_0$ to a point in $S_i\times a$ in the obvious way to a
path from $m_0$ to a point in $S_i\times (a+\delta)$.
 In particular  we get an induced map
\[ f_\delta:H_i^\b(S_j\times a\subset M;R^k)\to H_i^\b(S_j\times (a+\delta)\subset M;R^k).\]
 Using our choices of paths and using Lemma \ref{lem:twiiso}  we get a
map
\[ H_i^\b(S_j;R^k)= H_i^\b(S_j\times 0;R^k)\xrightarrow{\cong} H_i^\b(S_j\times 0\subset M;R^k)
\to H_i^\b(N'\subset M;R^k) \xrightarrow{\cong} H_i^\b(N';R^k)\] which we denote by $\i_-$.
Similarly we get a map
\[ \ba{rcl} H_i^\b(S_j;R^k)&=& H_i^\b(S_j\times 0;R^k)\\
&\xrightarrow{\cong}& H_i^\b(S_j\times 0\subset M;R^k)
\xrightarrow{f_{w_j-1+\eps}}H_i^\b(S_j\times (w_j-1+\eps)\subset M;R^k)\\
&\to& H_i^\b(N' \subset M;R^k) \xrightarrow{\cong} H_i^\b(N';R^k)\ea \] which we denote by
$\i_+$.

For the remainder of this section we pick  a representation $\a:\pi_1(M,m_0)\to \glfk$.  With
our conventions and choices we can now formulate the following crucial lemma which provides the
link between the twisted homology of $S_1,\dots,S_l$ and the homology of $M$.

\begin{proposition}\label{prop:exact}
There exists a long exact  sequence
\[ \dots \to  \bigoplus_{j=1}^l H_i^{\a}(S_j;\fk)\otimes_{\F} \ft \xrightarrow{
\bigoplus\limits_{j=1}^l \i_--\i_+t^{w_j}} H_i^{\a}(N';\fk)\otimes_{\F} \ft
   \to H_i^{\a\otimes \phi}(M;\fkt)
  \to \dots \]
\end{proposition}

\begin{proof}
We have the following Mayer--Vietoris sequence of twisted subspace homology (where we write
$V=\fkt$):
\[
\to \ba{c} H_i^{\a\otimes \phi}(S\times \eps\subset M;V)\\\oplus \\H_i^{\a\otimes \phi}(S\times
0\subset M;V)\ea \hspace{-0.1cm} \xrightarrow{\tiny{\big( \hspace{-0.1cm}\ba{cc} \i&\hspace{-0.2cm}
\i\\-\i&\hspace{-0.2cm}-\i\ea \hspace{-0.1cm}\big)}} \hspace{-0.3cm} \ba{c} H_i^{\a\otimes
\phi}(N\subset M;V)\\\oplus
\\H_i^{\a\otimes \phi}(S\times [0,\eps]\subset M;V)\ea \hspace{-0.2cm} \xrightarrow{\tiny{( \ba{cc}\i & \i\ea)}}
H_i^{\a\otimes \phi}(M;V) \to\] where $\i$ stands for the maps induced by the various injections.
Now consider the following commutative diagram of sequences
\[ \ba{cccccccccccccccccc}
\ba{c} H_i^{\a\otimes \phi}(S\times \eps\subset M;V)\\\oplus \\H_i^{\a\otimes \phi}(S\times
0\subset M;V)\ea
&\hspace{-0.2cm}\xrightarrow{\tiny{\bp \i &\i\\
-\i&-\i\ep}}&\hspace{-0.4cm} \ba{c} H_i^{\a\otimes \phi}(N\subset M;V)\\\oplus \\H_i^{\a\otimes
\phi}(S\times [0,\eps]\subset M;V)\ea\hspace{-0.4cm}& \hspace{-0.2cm}\xrightarrow{\tiny{(
\ba{cc}\i & \i\ea)}}& \hspace{-0.2cm}H_i^{\a\otimes \phi}(M;V)\\
 \Big \uparrow
{\Big(\hspace{-0.1cm} \ba{c} -f_\eps\\\id \ea\hspace{-0.1cm} \Big)}
&&\Big \uparrow {\Big(\hspace{-0.1cm} \ba{c} \id\\0\ea\hspace{-0.1cm} \Big)} &&\Big \uparrow&\\
 H_i^{\a\otimes \phi}(S\times 0\subset M;V) \hspace{-0.4cm}&\hspace{-0.4cm}
\xrightarrow{ \i-\i \circ f_\eps}& \hspace{-0.4cm}H_i^{\a\otimes \phi}(N\subset M;V)
&\xrightarrow{\i}& \hspace{-0.42cm} H_i^{\a\otimes \phi}(M;V).\ea \] Note that given $a\in
H_i^{\a\otimes \phi}(S\times \eps\subset M;\fkt)$ and $b\in H_i^{\a\otimes \phi}(S\times 0\subset
M;\fkt) $ we have $\i(a)+\i(b)=0\in H_i^{\a\otimes \phi}(S\times [0,\eps]\subset M;\fkt)$ if and
only if $a=-f_{\eps}(b)$. From this it now follows easily that the bottom sequence is also exact.

Now note that we have canonical isomorphisms
\[ \ba{rcccl}
H_i^{\a\otimes \phi}(N\subset M;V)&\hspace{-0.2cm}\cong\hspace{-0.2cm}&H_i^{\a\otimes
\phi}(N'\subset M;V)&\hspace{-0.2cm}\oplus\hspace{-0.2cm}&\bigoplus\limits_{j=1}^l
\bigoplus\limits_{s=1}^{w_j-1} H_i^{\a\otimes \phi}(S_j\times
[s-1+\eps,s]\subset M;V)\\
H_i^{\a\otimes \phi}(S\subset M;V)&\hspace{-0.2cm}\cong\hspace{-0.2cm}&\bigoplus\limits_{j=1}^l
H_i^{\a\otimes \phi}(S_j\times 0\subset
M;V)&\hspace{-0.2cm}\oplus\hspace{-0.2cm}&\bigoplus\limits_{j=1}^l
\bigoplus\limits_{s=1}^{w_j-1} H_i^{\a\otimes \phi}(S_j\times s\subset M;V).\ea \] It is now
easy to see, using the same arguments as above,
  that the following sequence is
exact as well:
\[
\to \bigoplus_{j=1}^l H_i^{\a\otimes \phi}(S_j\times 0\subset M;V) \xrightarrow{
\bigoplus\limits_{j=1}^l \i-\i \circ f_{w_j-1+\eps}} H_i^{\a\otimes \phi}(N'\subset M;V)
\xrightarrow{\i} H_i^{\a \otimes \phi}(M;V) \to \dots.\]

Now note that $\phi$ vanishes on  $H_1(N')$ and on every $H_1(S_j)$. Indeed, every curve in
$S_j$ can be pushed off into $N'$, where $\phi$ vanishes. We therefore get {\emph{canonical}}
isomorphisms $H_i^{\a\otimes \phi}(N';\fkt)\cong H_i^\a(N';\fk)\otimes_\F \ft$ and
$H_i^{\a\otimes \phi}(S_j;\fkt)\cong H_i^\a(S_j;\fk)\otimes_\F \ft$. We are done once we prove
the following claim.
\begin{claim}
The diagram
\[ \ba{cccccccccccccccccc}
 H_i^{\a\otimes \phi}(S_j\times 0\subset M;\fkt) &\xrightarrow{ \i-\i \circ f_{w_j-1+\eps}}& H_i^{\a\otimes \phi}(N'\subset M;\fkt)\\
\Big \downarrow {\cong}&& \Big \downarrow {\cong}&&\\
 H_i^{\a\otimes \phi}(S_j\times 0;\fkt) &\longrightarrow& H_i^{\a\otimes \phi}(N';\fkt)\\
\Big \downarrow {\cong}&& \Big \downarrow {\cong}&&\\
H_i^\a(S_j;\fk)\otimes_{\F}\ft &\xrightarrow{ \i_--\i_+t^{w_j}}& H_i^\a(N';\fk)\otimes_{\F}\ft .\ea
\] commutes. Here the top vertical maps are given by the isomorphisms induced from the choice of
paths, and the bottom  isomorphisms are the canonical isomorphisms mentioned above.
\end{claim}

In order to prove the claim first recall that $\ti{M}$ can be viewed as the homotopy classes of
paths in $M$ starting at $m_0$. Since $\phi$ vanishes on $N'$ we can define $\phi:\ti{N'}\to \Z$
(recall that $\ti{N'}\subset \ti{M}$) by sending $q\in \ti{N'}$ (represented by a path which we
also call $q$) to $\phi$ of the closed path given by juxtaposing $q$ with a path in $N'$ from
the end point of $q$ to $m_0$. This is a well--defined surjective map and we can decompose
$\ti{N'}=\cup_{r\in \Z} \ti{N'}_r$ where $\ti{N'}_r=\phi^{-1}(r)$. Now note that the isomorphism
from Lemma \ref{lem:twiiso} (applied to the constant path which connects the base point $m_0$ of
$N'$ with the base point $m_0$ of $M$) gives an isomorphism
\[ \ba{rcl}
H_i((\oplus_{r} C_*(\ti{N'}_r))\otimes_{\Z[\pi_1(M,m_0)]}\fkt)&=& H_i^{\a\otimes \phi}(N'\subset
M;\fkt)\\
& \cong &H_i^{\a\otimes \phi}(N';\fkt)\\
&\cong &H_i^\a(N';\fk)\otimes_{\f}\ft\ea \] where an element represented by $\s_r\otimes vt^l$ with
$\s_r\in C_*(\ti{N'}_r)$ and $v\in \fk$ gets sent to an element of the form $\s\otimes t^{r+l}$
where $\s\in H_i^\a(N';\fk)$.

Similarly we can decompose $\widetilde{S_j\times 0}=\cup_{r\in \Z} \widetilde{(S_j\times 0)}_r$
(using paths in $N'$ again), and the same arguments apply.

Now note that the inclusion $\widetilde{S_j\times 0}\to \widetilde{N'}$ clearly sends
$\widetilde{(S_j\times 0)}_r$ into $\widetilde{N'}_r$. On the other hand the map
$f_{w_j-1+\eps}:\widetilde{S_j\times 0}\to\widetilde{N'}$ sends a point in $\widetilde{S_j\
\times 0}$ represented by a path $\g$ to the the point represented by the extension of the path
$\g$ through $S_j\times [0,w_j-1+\eps]$. Closing it up in $N'$ we get a path whose intersection
number with $S$ is increased by $w_j$. This shows that $f_{w_j-1+\eps}$ sends
 $\widetilde{(S_j\times 0)}_r$ into $\widetilde{N'}_{r+w_j}$. The claim now follows easily
from the above observations.
\end{proof}

%Note that because of base point issues this induced
%homomorphism is only defined up to conjugacy. But the homology
%groups $H_*^{\a}(X;\F^k)$ are isomorphic, and their dimensions
%over $\F$ are well-defined.
%Homology with local coefficients depends on a choice of base points, but for path connected spaces
%different choices of base points give isomorphic homology groups. For non--connected spaces it is
%best to think of local coefficients as a locally constant sheaf, so that the restriction to any
%subspace (connected or not) gives well--defined homology groups. We will therefore suppress base
%points and the choice of paths connecting base points in our notation.

In the following we write $b_n^\a(S) := \dim_\F(H_n^\a(S;\fk))=\sum_{i=1}^l
w_i\dim_\F(H_n^\a(S_i;\fk))$.

\begin{proposition} \label{propbis1}
We have
\[ b_1^{\a}(S) \geq  \dimm_{\F}\left(\tor_{\ft}\left(\twihomphi\right)\right).\]
In particular if $\twialex \ne 0$, then $b_1^\a(S) \ge
\degtwialex$.
\end{proposition}

The proof is a variation on a standard argument.

\begin{proof}
Consider the exact sequence from Proposition \ref{prop:exact}. Note that
\[ F:= \ker\{\oplus_{j=1}^l
H_0^{\a}(S_j;\F^k)\otimes_\F \ft\to H_0^{\a}(N';\F^k)\otimes_\F \ft\}\subset \oplus_{j=1}^l
H_0^{\a}(S_j;\F^k)\otimes_\F \ft \] is a (possibly trivial) free $\ft$--module. Consider the
exact sequence
\[ \bigoplus\limits_{j=1}^l H_1^{\a}(S_j;\fk)\otimes_\F \ft \xrightarrow{ \bigoplus\limits_{j=1}^l
\i_--\i_+t^{w_j}} H_1^{\a}(N';\fk)\otimes_\F \ft \to H_1^{\a\otimes
\phi}(M;\fkt)\xrightarrow{\partial} F \to 0.
\] Since $\ft$ is a PID the map $\partial$ splits, i.e., $ H_1^{\a\otimes \phi}(M;\fkt)\cong
\ker(\partial)\oplus F$. In particular
\[ \dim_{\F}\left(\tor_{\ft}(\twihomphi)\right)=
\dim_{\F}\left(\tor_{\ft}(\ker(\partial))\right).\]
Using appropriate bases the map
$\bigoplus\limits_{j=1}^l \i_--\i_+t^{w_j} $, which represents the module $\ker(\partial)$,  is
presented by a matrix of the form
\[ \bp A_1t^{w_1}+B_1&\dots & A_lt^{w_l}+B_l\ep \]
where $A_j,B_j, j=1,\dots,l$ are matrices over $\F$ of size
$\dim_{\F}\left(H_1^{\a}(N';\fk)\right)\times \dim_{\F}\left(H_1^{\a}(S_j;\fk)\right)$. The
proposition now follows easily from combining Lemma \ref{lemmafinite} with the following claim.

\begin{claim}
 Let $H$ be a $\ft$--module with a
presentation matrix of the form
\[ C=\bp A_1t^{w_1}+B_1&\dots & A_lt^{w_l}+B_l\ep \]
where $A_j,B_j$ are matrices over $\F$ of size $p\times q_j$. Then
 $\dimm_{\F}(\tor_{\ft}(H))\leq \sum_{j=1}^l q_jw_j$.
\end{claim}

For the proof of the claim let $q=\sum_{j=1}^l q_j$. Using row and column operations over the
PID $\ft$ we can transform $C$ into a matrix of the form
\[ \bp f_1(t) &0&\dots&0&0 \\
 0 &f_2(t)&\dots&0&0 \\
0 &0&\ddots&0&0 \\
 0 &0&\dots&f_r(t)&0\\
 0&0&\dots&0&(0)_{p-r \times q-r} \ep \]
for some $f_i(t)\in \ft\sm \{0\}$. Clearly $\dim_{\F}(\tor_{\ft}(H))=\sum_{i=1}^r \deg(f_i(t))$.
Since row and column operations do not change the ideals of $\ft$ generated by minors (cf.
\cite[p.~101]{CF77}), and since any $k\times k$ minor of $At+B$ has degree at most $\sum_{j=1}^l
q_jw_j$, it follows that $\sum_{i=1}^l \deg(f_i(t))\leq \sum_{j=1}^l q_jw_j$. This concludes the
proof of the claim.
\end{proof}

\begin{proposition} \label{prop:connected}
If $\twialex \ne 0$ then either  $S$ is connected or $b_0^\a(S_j)=0$ for $j=1,\dots,l$.
\end{proposition}

The following proof is partly inspired by ideas of Turaev \cite{Tu02b}.

\begin{proof}
 Consider the exact sequence from Proposition
\ref{prop:exact}:
\[ \ba{rl} \to & \twihomphi \\
\to \bigoplus\limits_{j=1}^l H_0^\a({S_j};\fk)\otimes_\F \ft\xrightarrow{ \bigoplus\limits_{j=1}^l
\i_--\i_+t^{w_j}} H_0^\a(N';\fk)\otimes_\F \ft\to& \twihomphizero \to 0.\ea \] From $\twialex \ne
0$ it follows that $\twihomphi$ is $\ft$--torsion. By Lemma \ref{lem:delta03} $\twihomphizero$ is
$\ft$--torsion. The exact sequence shows that the ranks of the free $\ft$--modules $\oplus_{j=1}^l
H_0^\a({S_j};\fk)\otimes_\F \ft$ and $H_0^\a(N';\fk)\otimes_\F \ft$ are equal, and hence
%If we now consider the above exact sequence with
%$\f(t)$--coefficients it follows that
\be \label{h0equal} \oplus_{j=1}^l H_0^\a({S_j};\fk)\cong H_0^\a(N';\fk). \ee

%\noindent Since we can arrange $w_i$ parallel copies of $S_i$ inside $\nu(S_i)$ in $M$, we see that
%$N \cong (M\sm \nu |\hat{S}|) \, \cup \, \bigcup\limits_{i=1}^l \bigcup\limits_{j=1}^{w_i-1}
%S_i\times [-1,1]$. Therefore we have the following isomorphisms
%\be \label{equnisom} \ba{rccccccc}
%H_0^\a({S};\fk)&\cong & \bigoplus\limits_{i=1}^l  H_0^\a(S_i;\fk) &\oplus&  \bigoplus\limits_{i=1}^l  H_0^\a(S_i;\fk)^{w_i-1}\\
% H_0^\a(N;\fk)&\cong&H_0^\a(M\sm \nu |\hat{S}|;\fk)&\oplus&  \bigoplus\limits_{i=1}^l  H_0^\a(S_i;\fk)^{w_i-1}\ea \ee
%where $H_0^\a(S_i;\fk)^{w_i-1} := \bigoplus\limits^{w_i-1} H_0^\a(S_i;\fk)$.
 Note that the representations $ \pi_1(S_j,s_j)\to \pi_1(M,m_0)\xrightarrow{\a} \glfk$ induced by our chosen
paths factor through $\pi_1(N',m_0)$. Therefore
\be \label{inequ} b_0^\a(S_j)\geq b_0^\a(N'), j=1,\dots,l
\end{equation} by Lemma \ref{lemmah0surjects}.

First consider the case $b_0^\a(N')=0$. In that case it follows from isomorphism (\ref{h0equal})
 that $b_0^\a(S_j)=0$ for
all $j=1,\dots,l$.
%If $S_i$ is not closed, clearly $b_2^\a(S_i)=0$.
%If $S_i$ is closed, then by Poincar\'{e} duality (Lemma \ref{lemmaduality}) $b_2^\a(S_i) =
%b_0^\a(S_i)=0$.

Now assume  that $b_0^\a(N')>0$. It follows immediately from the isomorphism in (\ref{h0equal})
%and (\ref{equnisom})
and from the inequality (\ref{inequ}) that $l=1$. But since $\phi$ is
primitive it also follows that $w_1=1$, i.e., $S$ is connected.
\end{proof}

\begin{proposition} \label{prop:b0} Assume that $\twialex \ne 0$,
then
\[ b_0^{\a}(S)= \degtwialexzero.\]
\end{proposition}

\begin{proof}  First assume that $b_0^{\a}(S_j)=0$ for every
component $j=1,\dots,l$. Then $H_0^\a(N';\fk) = 0$ by the isomorphism (\ref{h0equal}). This
implies that $H_0^{\a\otimes \phi}(M;\fkt)=0$ from the exact sequence in Proposition
\ref{prop:exact}, hence $\twialexzero=1$.

Now assume that $b_0^\a(S_j) \ne 0$ for some $j$. By Proposition \ref{prop:connected} $S$ is
connected and hence $N'=N$.  It follows from Lemma \ref{lemmah0surjects} that the maps $\i_+,\i_-:
H_0^\a({S};\fk)\to H_0^\a(N';\fk)$ are surjective. Since $H_0^\a({S};\fk) \cong H_0^\a(N';\fk)$ by
isomorphism (\ref{h0equal}) it follows that $\i_+$ and $\i_-$ induce isomorphisms on
$H_0^\a({S};\fk)$.

Let $b:=b_0^\a(S)=b_0^\a(N')$. Picking appropriate bases for $H_0^\a({S};\fk)$ and $H_0^\a(N';\fk)$
the sequence from Proposition \ref{prop:exact} becomes
 \[ \f^b \otimes_\F \ft \xrightarrow{t\cdot Id-J } \f^b\otimes_\F \ft \to H_0^{\a\otimes \phi}(M;\fkt)\to 0,\]
 where $J:\f^b\to \f^b$ is an isomorphism. It follows that
$H_0^{\a\otimes \phi}(M;\fkt)\cong \f^b\cong H_0^\a(S;\fk)$. The lemma now follows from Lemma
\ref{lemmafinite}.
\end{proof}

%\noindent We note that from Propositions \ref{prop:connected} and \ref{prop:b0} we immediately get
%the following useful corollary:
%\begin{corollary} \label{cor:conn}
%If $\twialexzero \ne 1$ and $\twialex\ne 0$, then there exists a
%Thurston norm minimizing surface which is connected.
%\end{corollary}

\begin{proposition} \label{prop:b2}
Assume that $\partial M$ is empty or consists of tori. If $\twialex \ne 0$, then
\[ b_2^{\a}(S)= \degtwialextwo.\]
\end{proposition}

\begin{proof}
 Let $I:=\{ i\in
\{1,\dots,l\}| S_i \mbox{ closed} \}$ and let $T=\cup_{i\in I}S_i$.  Clearly
$b_2^\a(S)=\sum_{i\in I} w_ib_2^\a(S_i)$. Note that we can write $\partial N'=T_+ \cup T_-\cup
W$ for some surface $W$ where $T_-=T\times 0$ and $T_+=\cup_{i\in I}S_i\times (w_i-1+\eps)$. It
follows from Lemma \ref{lem:delta03} and Proposition \ref{prop:delta2} that the long exact
sequence from Proposition \ref{prop:exact} becomes
\[ 0 \to \bigoplus_{i\in I} H_2^\a(S_i;\fk)\otimes_\F \ft \xrightarrow{
\bigoplus_{i\in I} \i_--\i_+t^{w_i}} H_2^\a(N';\fk)\otimes_\F \ft\to \twihomphitwo \to 0.
\] We need the following claim.
\begin{claim} The maps $\i_-,\i_+:H_2^\a(T;\fk)=\oplus_{i\in I} H_2^\a(S_i;\fk)\to H_2^\a(N';\fk)$ are isomorphisms.
\end{claim}

In order to give a proof of the claim we first tensor the above short exact sequence with
$\f(t)$. We see that $H_2^\a(T;\fk)$ and $H_2^\a(N';\fk)$ have the same dimension as $\F$-vector
spaces. It is therefore enough to show that $\i_-$ and $\i_+$ are injections. This is clearly
the case if $T=\emptyset$. So let us now assume that $T\ne \emptyset$.

Consider the short exact sequence
\[ H_3^\a(N',T_+;\fk)\to H_2^\a(T_+;\fk) \to H_2^\a(N';\fk). \]
Note that $\partial N'$ is the disjoint union of $T_+$ and $T_-\cup W$ since $T_+$ is closed. We
can therefore apply Poincar\'e duality. By Poincar\'e duality and by Lemma \ref{lemmadualitybig} in
Section \ref{sectionduality} we then have
\[ H_3^\a(N',T_+;\fk)\cong H^0_\a(N',T_-\cup W;\fk)\cong \hom_\f(H_0^{\at}(N',T_-\cup W;\fk),\f). \]
\noindent Here $\at$ is the adjoint representation of $\a$ which is defined in Section
\ref{sectionduality}. Since $T\ne \emptyset$ it follows from Lemma \ref{lemmah0surjects} that the
map $H_0^{\at}( T_-\cup W;\fk) \to H_0^{\at}(N';\fk)$ is surjective. It now follows from the long
exact homology sequence that $ H_0^{\at}(N',T_-\cup W;\fk)=0$. This shows that $\i_+$ is injective.
The proof for $\i_-$ is identical. This concludes the proof of the claim.

We now showed that $\twihomphitwo$ has a presentation matrix of the form $AT+B$ where $A,B$ are
invertible matrices over $\F$ and $D$ is a diagonal matrix with $b_2^\a(S_i)$ entries $t^{w_i}$
for any $i\in I$. Note that
\[ \det(AD+B)=\det(B)+\dots+\det(A)t^{\sum_{i\in I}w_ib_2^\a(S_i)}.\]
It follows that
\[ \dim(\twihomphitwo)=\deg(\det(AD+B))=\sum_{i\in I}w_ib_2^\a(S_i)=b_2^\a(S).\]
\end{proof}

%=================
Now we can conclude the  proof of Theorem \ref{mainthm}.

\begin{proof}[{\bf Proof of Theorem \ref{mainthm}}]
Without loss of generality we can assume that $\phi$ is primitive since the Thurston norm and the
degrees of twisted Alexander polynomials are homogeneous. Let $\hat{S} $ be the weighted surface
from Proposition \ref{propminus}.  By Lemma \ref{lem:euler} we have
\[ \ba{rcl} ||\phi||_T&=&\max\{0,b_1(S)- (b_0(S) +
 b_2(S))\}\\[1mm]
&\geq& b_1(S)- (b_0(S) + b_2(S))\\[1mm]
&=&\frac{1}{k}\left(b_1^\a(S)-\left(b_0^{\a}(S)+b_2^{\a}(S)\right)\right). \ea\]
 The theorem now follows immediately from Propositions
\ref{propbis1},  \ref{lem:delta03}, \ref{prop:b0}, \ref{prop:b2} and Lemma \ref{lemmareidemeister}.

%===========================================
\subsection{Proof of Theorem \ref{mainthm2}} \label{sectionfibered} \label{sectionfib}

\begin{proof}[Proof of Theorem
\ref{mainthm2}] Let $S$ be a fiber of the fiber bundle $M\to S^1$. Clearly $S$ is dual to $\phi \in
H^1(M)$ and it is well--known that $S$ is Thurston norm minimizing. Denote by $\hat{M}$ the
infinite cyclic cover of $M$ corresponding to $\phi$. Then an easy argument shows that
$H_i^{\a\otimes \phi}(M;\fkt)\cong H_i^\a(\hat{M};\fk)$ (cf. also \cite[Theorem 2.1]{KL99}). In
particular $H_i^{\a\otimes \phi}(M;\fkt)\cong H_i^\a(S;\fk)$.

 By assumption $S\ne D^2$ and $S\ne S^2$. Therefore by Lemmas
\ref{lem:euler} and \ref{lemmafinite} we get
\[ \ba{rcll} ||\phi||_T&=&\chi_-(S)&\\[1mm]
&=&b_1(S)-b_0(S)-b_2(S)&\\[1mm]
&=&\frac{1}{k} \left(b_1^\a(S)-b_0^\a(S)-b_2^\a(S)\right)&
%\mbox{ by Lemma \ref{lem:euler}}
\\[1mm]
&=&\frac{1}{k} \left(\dim_\f\big(H_1^{\a\otimes \phi}(M;\fkt)\big) -\dim_\f\big(H_0^{\a\otimes
\phi}(M;\fkt)\big)-\dim_\f\big(H_2^{\a\otimes \phi}(M;\fkt)\big)\right)&
%\mbox{ by Lemma \ref{lemmah0m}(2)}
\\[1mm]
&=&\frac{1}{k} \big(\degtwialex-\degtwialexzero-\degtwialextwo \big)&\\[1mm]
&=&\deg(\tau(M,\phi,\a)).
%\mbox{by Lemma \ref{lemmafinite}}.
\ea \]
\end{proof}

Since $||\phi||_T$ might be unknown for a given example the following corollary to Theorem
\ref{mainthm2} gives sometimes a more practical fibering obstruction.

\begin{corollary} \label{corfibered}
 Let $M$ be a 3--manifold and $\phi\in H^1(M)$ primitive such
that $(M,\phi)$ fibers over $S^1$ and such that $M\ne S^1\times D^2, M\ne S^1\times S^2$. Let $\F$
and $\F'$ be fields. Consider the untwisted Alexander polynomial $\alex \in \ft$. For any
representation $\a:\pi_1(M)\to \gl(\f',k)$ we have
 \[ \deg(\alex)-\oneplusbm=\frac{1}{k}\big(\degtwialex-
 \degtwialexzero-\degtwialextwo \big).\]
\end{corollary}

\begin{proof}
The corollary follows immediately from applying Theorem \ref{mainthm2} to the trivial
representation $\pi_1(M)\to \gl(\F,1)$ and to the representation $\a$.
\end{proof}

%
% \[ \ba{rcll} ||\phi||_T&=&\max\{0,b_1(S)- (b_0(S) +
% b_2(S))\}&\\[1mm]
%&\geq& b_1(S)- (b_0(S) + b_2(S))&\\[1mm]
%&=&1/k\left(b_1^\a(S)-\left(b_0^{\a}(S)+b_2^{\a}(S)\right)\right)&\mbox{
%by Lemma \ref{lemmab1}}\\[1mm]
%&=&1/k\left(b_1^\a(S)-\left(b_0^{\a}(S)+\ti{b}_3^\a(M)b_0^\a(S)\right)\right)
%&\mbox{ by Lemma \ref{prop:delta2}}(2)\\[1mm]
%&=&1/k\left(b_1^\a(S)-\left(\degtwialexzero+\ti{b}_3^\a(M)\degtwialexzero\right)\right)
%&\mbox{ by Lemma \ref{prop:b0}}.
%
%\ea\] The theorem now follows from Proposition \ref{propbis1} and
%Lemmas \ref{prop:delta2}(1) and \ref{lemmab3zero}.
\end{proof}

\section{The case of vanishing Alexander polynomials}\label{sectiontorsion}
Let $L$ be a boundary link (for example a split link). It is well--known that the untwisted
onevariable and multivariable Alexander polynomials of $L$ vanish (cf. \cite{Hi02}). Similarly one
can see that in fact most of the twisted onevariable and multivariable  Alexander polynomials
vanish as well. (See \cite{FK05} for the definition of twisted multivariable Alexander
polynomials.) Therefore Theorem \ref{mainthm} can in most cases not be applied to get lower bounds
on the Thurston norm.

It follows clearly from Propositions \ref{propbis1} and \ref{prop:b2} that the condition $\twialex
\ne 0$ is only needed to ensure that there exists a surface $S$ dual to $\phi$ with
$b_0^\a(S)=\degtwialexzero$ and $b_2^\a(S)=\degtwialextwo$. The following theorem can often be
applied in the case of link complements.

\begin{theorem} \label{thmthurstonlink} Let $M$ be a 3--manifold such that $H^1(M)\xrightarrow{i^*}
H^1(\partial M)$ is an injection where $i^*$ is the inclusion--induced homomorphism.  Let $N$ be a
torus component of $\partial M$, $\phi\in H^1(N) \cap \im{(i^*)}$ primitive, and $\a:\pi_1(M)\to
\glfk$ a representation.  Then
%there exists a connected surface $S\subset M$ such that $\partial S$ is dual to $\phi$ and
%\[ ||(i^*)^{-1}(\phi)||_{T,M} \geq  \frac{1}{k}\deg(\ti{\Delta}^\a_\phi(t))-\frac{1}{k}b_0^{\a}(S). \]
%In particular
\[ ||(i^*)^{-1}(\phi)||_{T} \geq  \frac{1}{k}\deg(\ti{\Delta}^\a_1(t))-1. \]
\end{theorem}

It is not hard to show that we can find a Thurston norm minimizing surface dual to
$(i^*)^{-1}(\phi)$ which is connected and has boundary (cf. e.g. \cite[Corollary 10.4]{Ha05} or
Turaev \cite[p.~14]{Tu02b}). The theorem now follows from the proof of Theorem \ref{mainthm}.

%\begin{proof}
%Let us consider the following commutative diagram
%\[ \ba{ccc} H^1(M)&\hookrightarrow &H^1(\partial M)\\
%\downarrow&&\downarrow \\
%H_2(M,\partial M)&\hookrightarrow &H_1(\partial M),\ea \] where the vertical maps are given by
%Poincar\'e duality. An embedded surface $S$ in $M$ is dual to $(i^*)^{-1}(\phi)$ if and only if
%$\partial S$ is dual to $\phi$. It follows that
%\[
%||(i^*)^{-1}(\phi)||_{T,M}=\min\{\chi_-(S)| S \mbox{ properly
%embedded}, \partial S \mbox{ Poincar\'e dual to }\phi \in
%H^1(\partial M)\}.
%\]
%Let $c\subset N$ be a simple closed curve Poincar\'e dual to $\phi$. Using that $N$ is a torus
%we can easily show that
%\[  ||(i^*)^{-1}(\phi)||_{T,M}=\min\{\chi_-(S)|S \mbox{ properly embedded and } \partial S=c \}.\]
%Let $S$ be a (possibly disconnected) Thurston norm minimizing surface  such that $\partial S =
%c$. By throwing away the components of $S$ which do not contain $c$, we can assume that $S$ is
%connected. Hence  $b_0(S)=1$ and $b_2(S)=0$. In particular $b_0^\a(S)\leq k$ and $b_2^\a(S)=0$.
%Therefore using Proposition \ref{propbis1} and Lemmas \ref{lemmafinite} and \ref{lem:euler} we
%obtain that
%\[ \ba{rcl} ||\phi||_T&\geq &b_1(S)-b_0(S)-b_2(S)\\
%&=&\frac{1}{k}\left(b_1^\a(S)-b_0^\a(S)-b_2^\a(S)\right)\\
%&\geq &\frac{1}{k}\deg(\ti{\Delta}^\a_1(t))-1.\ea \]
%\end{proof}

The main application is to study the Thurston norm of the complement of a link $L = L_1\cup \dots
\cup L_m\subset S^3$. In this case we can take $\phi$ to be dual to the meridian of the $i^{th}$
component $L_i$. Then it follows from the proof of Theorem \ref{thmthurstonlink} and a standard
argument that $ ||(i^*)^{-1}(\phi)||_{T}=2 \, \gen(L_i)-1$, where $\gen(L_i)$ denotes the minimal
genus of a surface in $X(L)$ bounding a parallel copy of $L_i$. Similar results were obtained by
Turaev \cite[p.~14]{Tu02b} and Harvey \cite[Corollary 10.4]{Ha05}.

%Using work of Freedman and He \cite{FH91} Cantarella, Kusner and Sullivan \cite[Corollary
%22]{CKS02} showed that the Thurston norm can be used to give lower bounds on the ropelength of a
%link component. They formulated a conjecture for a certain link. This conjecture was proved by
%Harvey \cite[Section 8]{Ha05} using higher--order Alexander polynomials. Using one--dimensional
%representations together with Theorem \ref{thmthurstonlink} and Lemma \ref{lemmaeven} we can
%reprove this conjecture.

The following observation will show that in more complicated cases there is no immediate way to
determine $b_0(S)$: if $L=L_1\cup L_2$ is a split oriented link, and $\phi:H_1(X(L))\to \Z$ given
by sending the meridians to 1, then a Thurston norm minimizing surface $S$ dual to $\phi$ is easily
seen to be  the disjoint union of the Seifert surfaces of $L_1$ and $L_2$.
%This follows immediately
%from the proof that the genus is additive, i.e., $\gen(L_1\#L_2)=\gen(L_1)+\gen(L_2)$ (cf.
%\cite[p.~18]{Lic97}).
%In particular $b_0(S)=2$.
On the other hand if $L_1$ and $L_2$ are parallel copies of a knot with opposite orientations and
$\phi:H_1(X(L))\to \Z$ is again given by sending the meridians to 1, then the annulus $S$ between
$L_1$ and $L_2$ is dual to $\phi$ with Euler characteristic zero.
%In particular it is connected, hence $b_0(S)=1$.
Summarizing, we have two situations in which  the first twisted Alexander polynomials vanish (in
fact $H_1(X(L);\qt)$ has rank one), $\phi$ is of the same type, but $b_0(S)$ differs.
%\footnote{I
%personally think that this remark is reasonably interesting, so I have an inclination towards
%keeping it in the paper, but if you think it is fluffy (or appears fluffy), then we can erase it}
%\footnote{I thought about adding a remark on the ropelength example, but finally decided against
%it}

%===========================================
%\section{Applications} \label{sectionapp}

%===========================================
\section{Examples} \label{sectiondet1}\label{sectionexamples}

%We applied our results to many explicit situations. In all
%reasonable situations we found the correct Thurston norm bounds
%and we found whether a manifold fibers or not.

%===========================================
\subsection{Representations of 3--manifold groups} \label{sectionrep}
%The study of the representation theory of 3--manifold groups has a long history. For example after the introduction of Casson's invariant
%$SU(2)$--representations were studied intensively.

%Let $M$ be a 3--manifold. Assume we are given a presentation $\langle g_1,\dots,g_s|
%r_1,\dots,r_t\rangle$ for $\pi_1(M)$. Then finding a representation to $\glfk$ for some $k$ is easy
%in theory: it is enough to assign arbitrary elements in $\glfk$ to $g_1,\dots,g_s$ and check
%whether these satisfy the relations. Our experience shows that this is not an effective way of
%finding representations since $\gl(\F_p,k)$ has approximately $p^{k^2}$ elements, and therefore
%there are $\big(p^{k^2}\big)^s$ possible assignments of elements in $\gl(\F_p,k)$ to $s$
%generators.

In our applications we first find homomorphisms $\pi_1(M)\to S_k$,
and then find a representation of $S_k$. Here $S_k$ denotes the
permutation group of order $k$.
%In all our examples we take $G=S_k$ for some $k$. (Metabelian groups can also be useful.)
The first choice of a representation for $S_k$ that comes to mind is $S_k\to \gl(\F,k)$ where $S_k$
acts by permuting the coordinates. But $S_k$ leaves the subspace $\{ (v,v,\dots,v) | v\in
\F\}\subset \F^k$ invariant, hence this representation is `not completely non--trivial'.
%It follows from Maschke's theorem
%that if $\gcd(\char(\F),|S_k|)=1$, then $\a$ decomposes as a direct sum of the trivial
%representation and a representation $V_{k-1}$.
To avoid this we prefer to work with a slightly different representation of $S_k$. If
$\varphi:\pi_1(M)\to S_k$ is a homomorphism then we consider
 \[
\a(\varphi):\pi_1(M)\xrightarrow{\varphi} S_k\to \gl(V_{k-1}(\F)),\]
 where
 \[ V_l(\F):=\{
(v_1,\dots,v_{l+1})\in \F^{l+1} | \sum_{i=1}^{l+1} v_i =0 \}.\]
 Clearly $\dim_\f(V_l(\f))=l$ and $S_{l+1}$
acts on it by permutation.
%Since $\a(\varphi)$ is a
%subrepresentation of a unitary representation, $\a(\varphi)$ is
%also unitary.
%These representations are easy to find and
%remarkably useful for our purposes.
%\\

We point out that the fundamental groups of 3--manifolds for which
the geometrization conjecture holds are residually finite (cf.
\cite{Th82} and \cite{He87}). In particular most (or perhaps all)
3--manifolds have many homomorphisms to finite groups, hence to
$S_k$'s.
%and in particular they have many
%interesting finite representations.

%We quickly explain why this approach is promising. Recall that irreducible manifolds with
%$b_1(M)\geq 1$ are Haken. Thurston \cite{Th82} (cf. also \cite{He87}) showed that the fundamental
%group of a Haken manifold is residually finite. Recall that a group $G$ is called \emph{residually
%finite} if for every $g\ne e \in G$ there exists a homomorphism $\a:G\to H$, $H$ a finite group,
%such that $\a(g)\ne e$. Furthermore the free product of residually finite groups is residually
%finite. This shows that any manifold we are interested in has many homomorphisms to finite groups.
%In fact the geometrization conjecture implies that all 3--manifold groups are residually finite.
%
%Note that every finite group $G$ is a subgroup of $S_{|G|}$. In particular the homomorphisms to
%$S_k, k\in \N$, contain all homomorphisms to all finite groups.

%===================================
\subsection{Knots with up to 12 crossings: genus bounds and fiberedness} \label{sectionexamples11}
\label{section11n34} In this section we show that the degrees of twisted Alexander polynomials
detect the genus of all knots with 12 crossings or less. Also we detect all non--fibered knots
with 12 crossings or less, some of which are new discoveries to our knowledge.
\\ \\
{\bf I. Knot genus.}
%Given a diagram for a knot one can find a
%Seifert surface using Seifert's algorithm (cf. \cite{Rol90}). This
%gives an upper bound on the genus of a knot. It turns out that for
%all knots with 10 or fewer crossings the (untwisted) Alexander
%norm determines the knot genus, that is, we have \[ 2 \,
%\gen(K)=\deg(\Delta_K(t)).\] This equality also holds for all
%alternating knots (cf. \cite{Cr59,Mu58a,Mu58b}). On the other hand
%it is known that
%\[  2 \, \gen(K)>\deg(\Delta_K(t))\] for many knots with more than 10 crossings. We will discuss
%all 11--crossing and all 12--crossing knots with this property in this section.
%The genera of knots with 12 crossing or less are known (cf. \cite{CL}). We will show that for
%these knots the correct genus is detected by Theorem \ref{mainthm}.
There are 36 knots with 12 crossings or less for which $\genus(K)>\frac{1}{2}\deg\Delta_K(t)$ (cf.
e.g. \cite{CL} or \cite{Sto}). The most famous and interesting examples are $K=11_{401}$ (the
Conway knot) and $11_{409}$ (the Kinoshita--Terasaka knot). Here we use the \emph{knotscape}
notation.

%Using geometric methods Gabai \cite{Ga84} showed  that the genus of the Conway knot is 3 and
%that the genus  of the Kinoshita--Terasaka knot is two. The computation of the genus for all
%11--crossing knots was done by Jacob Rasmussen, using a computer assisted computation of the
%Oszv\'ath--Szab\'o knot Floer homology (cf. also \cite{OS04a} and \cite{OS04b}).
%%, namely $K=11^n_{34}$ (the
%%Conway knot), $11^n_{42}$ (the Kinoshita--Terasaka knot),
%%$11^n_{45}, 11^n_{67},$ $ 11^n_{73}, 11^n_{97}, 11^n_{152}$.

\begin{figure}[h] \begin{center}
\begin{tabular}{cc}  \includegraphics[scale=0.27]{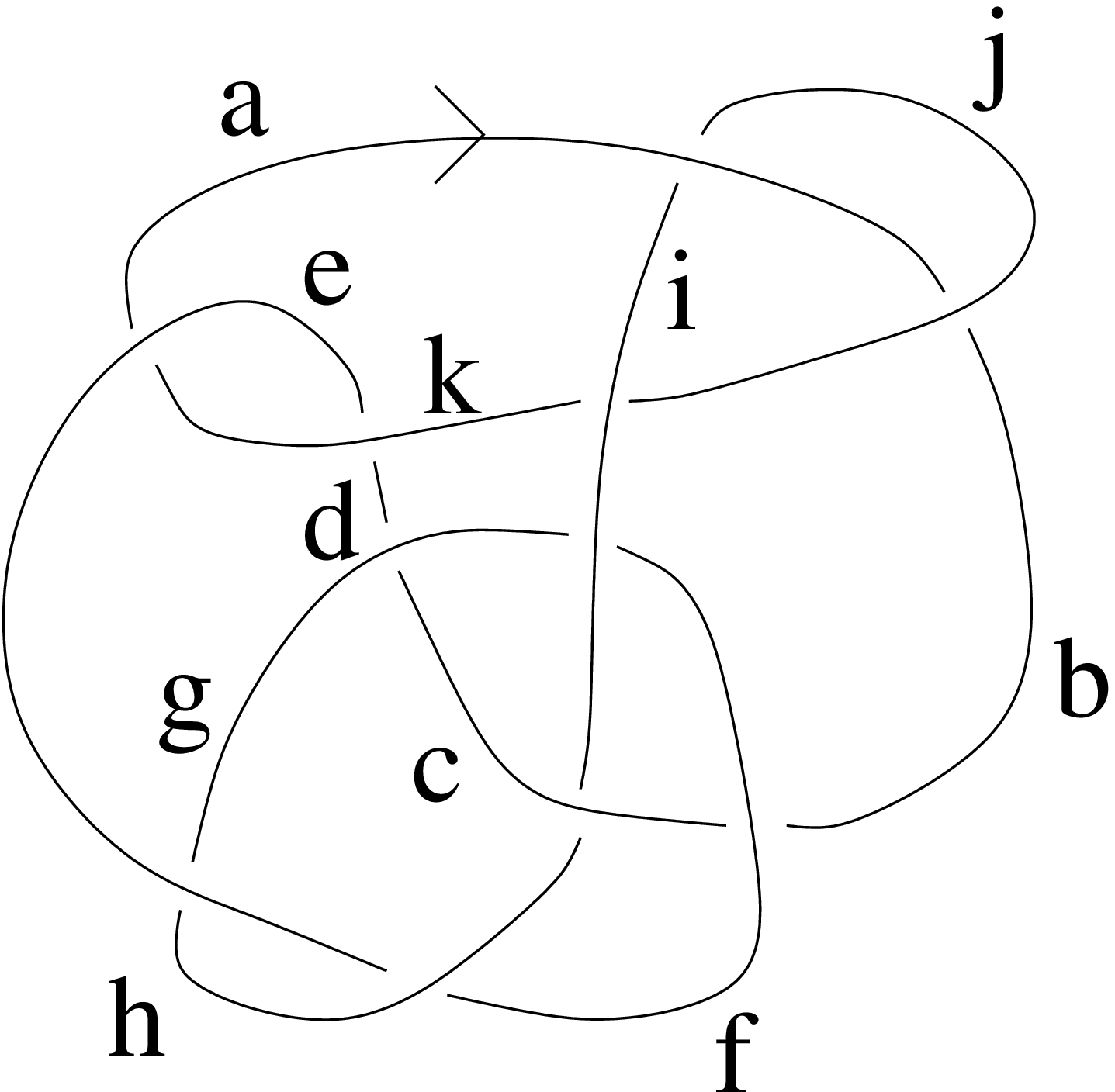}&\hspace{1cm}
\includegraphics[scale=0.28]{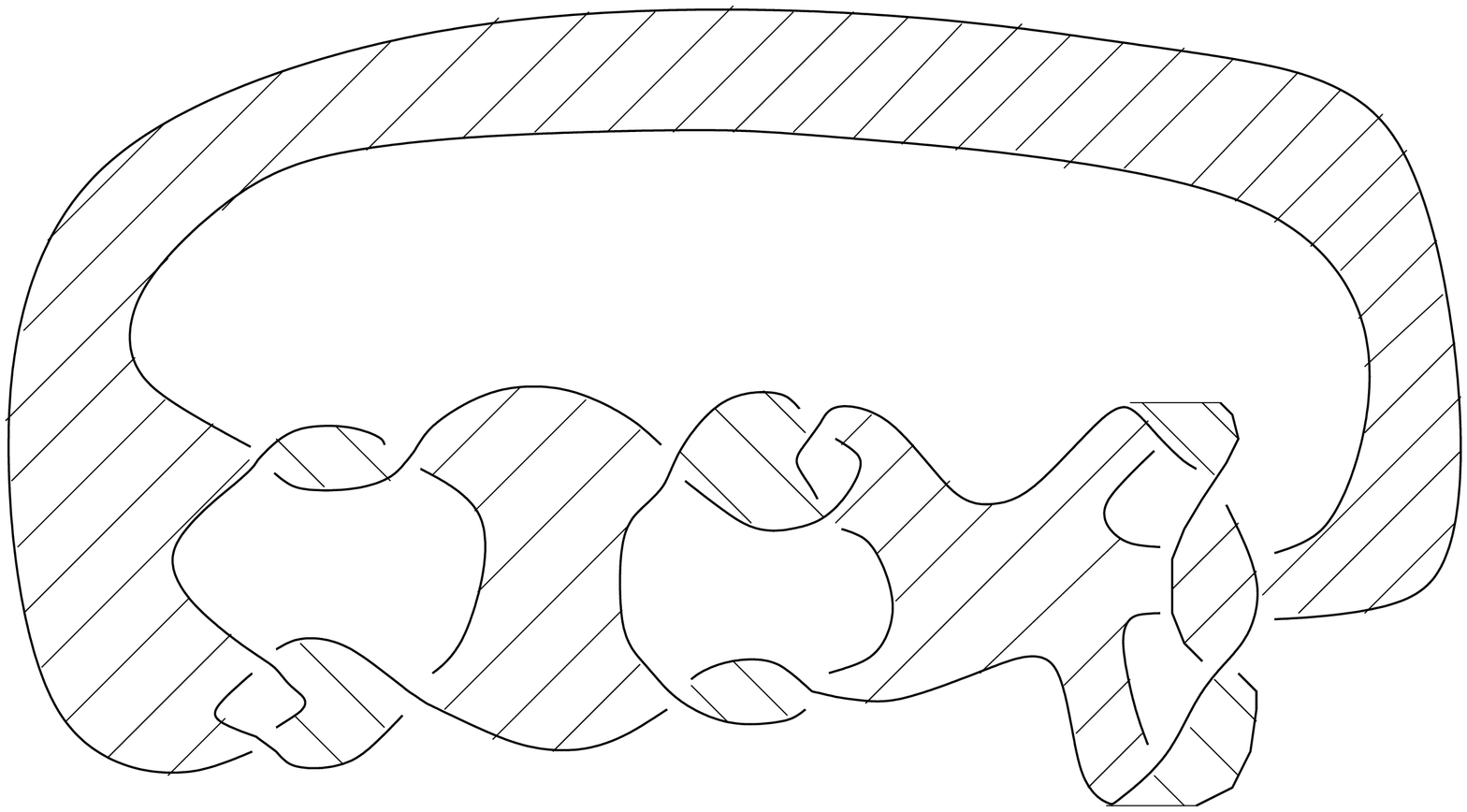} \end{tabular} \caption{The Conway knot
$11_{401}$ and a Seifert surface of genus 3 (from \cite{Ga84}).} \label{example11n34}
\end{center}
 \end{figure}

First, we consider the Conway knot $K=11_{401}$ whose diagram is given in Figure
\ref{example11n34}. The genus of the Conway knot is 3. This knot has Alexander polynomial one,
i.e., the degree of $\Delta_{K}(t)$ equals zero. Furthermore this implies that $\pi_1(X(K))^{(1)}$
is perfect, i.e., $\pi_1(X(K))^{(n)}=\pi_1(X(K))^{(1)}$ for any $n>1$. (For a group $G$, $G^{(n)}$
is defined inductively as follows: $G^{(0)} := G$ and $G^{(n+1)} := [G^{(n)},G^{(n)}]$.) Therefore
the genus bounds of Cochran \cite{Co04} and Harvey \cite{Ha05} vanish as well.
The fundamental group
$\pi_1(X(K))$ is generated by the meridians $a,b,\dots,k$ of the
segments in the knot diagram of Figure \ref{example11n34}. The
relations are \[ \ba{rclrclrclrclrcl}
a&=&jbj^{-1},&    b&=&fcf^{-1},&    c&=&g^{-1}dg,&    d&=&k^{-1}ek,&   \\
e&=&h^{-1}fh,&    f&=&igi^{-1},&    g&=&e^{-1}he,&    h&=&c^{-1}ic,\\
i&=&aja^{-1},&    j&=&iki^{-1},&    k&=&e^{-1}ae.& \ea \] Using the program \emph{KnotTwister}
\cite{Fr05}  we   found the homomorphism  $\varphi:\pi_1(X(K))\to S_5$ given by
 \[
\ba{rclrclrclrcl}
a&\mapsto&(142),& b&\mapsto&(451),& c&\mapsto&(451),& d&\mapsto&(453),\\
e&\mapsto&(453),& f&\mapsto&(351),& g&\mapsto&(351),& h&\mapsto&(431),\\
i&\mapsto&(351),& j&\mapsto&(352),& k&\mapsto&(321),&\ea \]
 where we use cycle notation.
%The generators of $\pi_1(X(K))$ are sent to the element in $S_5$
%given by the cycle with the corresponding capital letter.
We then consider $\a:=\a(\varphi):\pi_1(X(K))\xrightarrow{\varphi} S_5\to \gl(V_4(\F_{13}))$. Using
\emph{KnotTwister} we compute $\degtwialexzero=0$ and we compute the first twisted Alexander
polynomial to be
\[ \Delta_1^{\a}(t)= 1+6t+9t^2+12t^3+t^5+3t^6+t^7+3t^8+t^9
+12t^{11}+9t^{12}+6t^{13}+t^{14}\in \F_{13}[t^{\pm 1}]. \] Theorem \ref{mainthm} together with
Proposition  \ref{prop:delta2} says that if $\twialex\ne 0$, then
\[ \gen(K)
\geq \frac{1}{2k}\big(\degtwialex-\degtwialexzero\big)+\frac{1}{2}.
\]
Therefore in our case we get
\[ \gen(K)\geq \frac{1}{8}\cdot 14+\frac{1}{2}=\frac{18}{8}=2.25.\] Since $\gen(K)$ is an
integer we get $\genus(K) \geq 3$. Since there exists a Seifert surface of genus 3 for $K$ (cf.
Figure \ref{example11n34}) it follows that the genus of the Conway knot is indeed 3.

Second, let $K$ be the Kinoshita--Terasaka knot $K=11_{409}$. The genus of $K$ is 2. We found a map
$\varphi:\pi_1(X(K))\to S_5$ such that $\Delta_{1}^{\a(\varphi)}(t)\in \f_{13}[t^{\pm 1}]$ has
degree 12 and $\deg\big(\Delta_{0}^{\a(\varphi)}(t)\big)=0$. It follows from Theorem \ref{mainthm}
that $\genus(K)\geq \frac{1}{8}\cdot 12+\frac{1}{2}=2$.
%A Seifert surface of genus two
%is given in \cite{Ga84}.
Note that in this case our inequality becomes equality, hence `rounding up' is not necessary. Our
table below shows that this is surprisingly often the case. This fact is  of importance in
\cite{FK05} where we study the Thurston norm of link complements.

Table \ref{tab:compgenus12} gives all knots with 12 crossings or
less for which $\deg(\alexk)< 2 \, \gen(K)$. We obtained the list
of these knots from Alexander Stoimenow's knot page \cite{Sto}.
%One can also find the genus of all these knots on his knot page.
We compute twisted Alexander polynomials using \emph{KnotTwister} and 4--dimensional
representations of the form $\a(\varphi):\pi_1(X(K))\xrightarrow{\varphi} S_5\to
\gl(V_{4}(\F_{13}))$. Our genus bounds from Theorem \ref{mainthm} give (by rounding up if
necessary) the correct genus in each case.
 %All (no longer true for S_5 representations!!!!!!!!!!!!!!!)
%of the representations which give us the correct genus bounds have the property that
%$\deg\big(\Delta_{0}^{\a(\varphi)}(t)\big)=0$.

% changed the dimension to 5,
% degrees are now no longer in sync with genus bounds
\begin{table}[h]
\[\ba{|r|c|c|c|c|c|c|c|c|}
\hline
\mbox{Knotscape name}&11_{401}&11_{409}& 11_{412}& 11_{434}& 11_{440}& 11_{464}\\
\mbox{genus bound from $\alexk$}&0&0&2&1&2&1\\
%\mbox{dimension of $\a(\varphi)$}&4&4&4&4&2&4\\
%\mbox{degree of $\Delta_1^{\a(\varphi)}(t)$}&19&12&20&12&10&12\\
\mbox{genus bound from $\Delta_1^{\a}(t)$}&2.25&2.00&3.00&2.00&3.00&2.00\\
\hline
\mbox{Knotscape name}&11_{519}&12_{1311}&12_{1316}&12_{1319}&12_{1339}&12_{1344}\\
\mbox{genus bound from $\alexk$}&2&1&2&1&1&2\\
%\mbox{dimension of $\a(\varphi)$}&4&3&2 &2 &4&4\\
%\mbox{degree of $\Delta_1^{\a(\varphi)}(t)$}&20&9&10&10&12&20\\
%\mbox{upper bound on $b_0^\a(\varphi(\ker(\phi)))$}&0&0&0&0&0&0\\
\mbox{genus bound from $\Delta_1^{\a(\varphi)}(t)$}&3.00&2.00&2.50&3.00&2.00&3.00\\
\hline
\mbox{Knotscape name}&12_{1351}&12_{1375}&12_{1412}&12_{1417}&12_{1420}&12_{1509}\\
%\mbox{degree of  $\Delta_K(t)$}&4&2&2&4&4\\
\mbox{genus bound from $\alexk$}&2&2&1&1&2&2\\
%\mbox{dimension of $\a(\varphi)$}&2&4&4&4&4&2\\
%\mbox{degree of $\Delta_1^{\a(\varphi)}(t)$}&10&20&12&20&20&10\\
%\mbox{upper bound on $b_0^\a(\varphi(\ker(\phi)))$}&0&0&0&0&0&0\\
\mbox{genus bound from $\Delta_1^{\a(\varphi)}(t)$}&3.00&3.00&2.00&3.00&3.00&3.00\\
\hline
\mbox{Knotscape name}&12_{1519}&12_{1544}&12_{1545}&12_{1552}&12_{1555}&12_{1556}\\
%\mbox{degree of  $\Delta_K(t)$}&4&4&4&4&2&2\\
\mbox{genus bound from $\alexk$}&2&2&2&2&2&1\\
%\mbox{dimension of $\a(\varphi)$}&4&4&4&4&2&2\\
%\mbox{degree of $\Delta_1^{\a(\varphi)}(t)$}&20&20&20&20&10&6\\
%\mbox{upper bound on $b_0^\a(\varphi(\ker(\phi)))$}&0&0&0&0&0&0\\
\mbox{genus bound from $\Delta_1^{\a(\varphi)}(t)$}&3.00&3.00&3.00&3.00&3.00&2.00\\
\hline
\mbox{Knotscape name}&12_{1581}&12_{1601}&12_{1609}&12_{1699}&12_{1718}&12_{1745}\\
%\mbox{degree of  $\Delta_K(t)$}&0&2&2&0&2&2\\
\mbox{genus bound from $\alexk$}&1&0&1&1&0&1\\
%\mbox{dimension of $\a(\varphi)$}&4&5&4&4&4&4\\
%\mbox{degree of $\Delta_1^{\a(\varphi)}(t)$}&12&13&12&12&12&12\\
%\mbox{upper bound on $b_0^\a(\varphi(\ker(\phi)))$}&0&0&0&0&0&0\\
\mbox{genus bound from $\Delta_1^{\a(\varphi)}(t)$}&2.00&1.25&2.00&2.00&2.00&2.00\\
\hline
\mbox{Knotscape name}&12_{1807}&12_{1953}&12_{2038}&12_{2096}&12_{2100}&12_{2118}\\
%\mbox{degree of  $\Delta_K(t)$}&4&4&4&4&4\\
\mbox{genus bound from $\alexk$}&1&2&2&2&2&2\\
%\mbox{dimension of $\a(\varphi)$}&4&2&4&4&2&4\\
%\mbox{degree of $\Delta_1^{\a(\varphi)}(t)$}&12&10&20&20&10&20\\
%\mbox{upper bound on $b_0^\a(\varphi(\ker(\phi)))$}&0&0&0&0&0\\
\mbox{genus bound from $\Delta_1^{\a(\varphi)}(t)$}&2.00&3.00&3.00&3.00&3.00&3.00\\
\hline \ea \]
 \caption{Computation of degrees of twisted Alexander
 polynomials.} \label{tab:compgenus12}
 \end{table}
Using \emph{KnotTwister} it takes only a few seconds to find such
representations and to compute the twisted Alexander polynomial.\\

\begin{remark}
Let $K_1$ and $K_2$ be knots and assume there exists an epimorphism $\varphi:\pi_1(X(K_1))\to
\pi_1(X(K_2))$. Simon asked (cf. question 1.12 (b) on Kirby's problem list \cite{Kir97}) whether
this implies that $\genus(K_1)\geq \genus(K_2)$. Let $\a:\pi_1(X(K_2))\to \glfk$ be a
representation. By \cite{KSW05} $\Delta_{K_2,1}^{\a}(t)$ divides $\Delta_{K_1,1}^{\a\circ
\varphi}(t)$. Together with Lemma \ref{lemmah0surjects} this shows that the genus bounds from
Theorem \ref{mainthm} for $K_1$ are greater than or equal to the bounds for $K_2$. Thus Theorem
\ref{mainthm} (together with the observation that Theorem \ref{mainthm} often detects the
correct genus) suggests an affirmative answer to Simon's question. This should also be compared
to the results
 in \cite{Ha06}.
\end{remark}

\begin{remark} There are situations when for a given manifold the degree of
the twisted Alexander polynomial for some representation gives a worse bound for the Thurston norm
than the degree of the untwisted Alexander polynomial. This should be compared to the situation of
\cite{Co04,Ha06,Fr06}: Cochran's and Harvey's sequence of higher order Alexander polynomials gives
a never decreasing sequence of lower bounds on the Thurston norm.
\end{remark}

%==================================
{\bf II. Fiberedness.}
%It is known that a knot $K$ with 11 or fewer crossings is fibered if and
%only if $K$ satisfies
% \be \label{abelianfib} \alexk \mbox{ is monic and } \deg(\alexk)=2 \, \gen(K).\ee
% Hirasawa and Stoimenow had started a program to find all non--fibered 12--crossing knots. They
%showed that the knots $12_{1498}$, $12_{1502}$, $12_{1546}$ and $12_{1752}$ are not fibered even
%though they satisfy condition (\ref{abelianfib}). Now
Consider the knot $K=12_{1345}$. Its Alexander polynomial equals $\alexk=1-2t+3t^2-2t^3+t^4$ and
its genus equals two, therefore $K$ satisfies Neuwirth's condition (\ref{abelianfib}) in Section
\ref{section:introfib}. It follows from Corollary \ref{corfibered} that if $K$ were fibered,
then for any field $\F$ and any representation $\a:\pi_1(M)\to \gl(\f,k)$ the following would
hold:
 \[
\deg(\alexk)=\frac{1}{k}\big(\degtwialexk-\degtwialexzero\big)+1.\]
 We  found a representation
$\a:\pi_1(X(K))\to S_4$ such that for the canonical representation $\a:\pi_1(X(K))\to S_4 \to
\gl(\F_3,4)$ given by permuting the coordinates, we get $\deg(\Delta_1^\a(t))=7$ and $\degtwialexzero=1$. We compute
 \[
\frac{1}{4}\degtwialexk-\frac{1}{4}\degtwialexzero+1=\frac{10}{4}\ne
4 = \deg(\alexk).\] Hence $K$ is not fibered.

Similarly we found altogether 13 12--crossings knots which satisfy condition (\ref{abelianfib})
but which are not fibered; we list them in Table \ref{tab:nonfib}.
 %We used Corollary \ref{corfibered}
%as above. That is, we compared the degrees of untwisted Alexander polynomials with the degrees of
%twisted Alexander polynomials corresponding to some representation $\pi_1(S^3\sm K)\to S_k\to
%\gl(\F_p,k)$ where $S_k \to \gl(\F_p,k)$ is the canonical representation.
\begin{table}[h]
\[\ba{|r|c|c|c|c|c|c|c|c|}
\hline
\mbox{Knotscape name}&12_{1345}&12_{1498}&12_{1502}&12_{1546}&12_{1567}&12_{1670}&12_{1682}\\
\mbox{Order of permutation group $k$}&4&5&5&3&5&5&4\\
\mbox{Order $p$ of finite field}&3&2&11&2&3&2&3\\
\hline
\mbox{Knotscape name}&12_{1752}&12_{1771}&12_{1823}&12_{1938}&12_{2089}&12_{2103}&\\
\mbox{Order of permutation group $k$}&3&3&5&5&5&4&\\
\mbox{Order $p$ of finite field}&2&7&7&11&2&3&\\
\hline \ea \] \caption{Alexander polynomials of non--fibered knots}\label{tab:nonfib}
\end{table}
As we mentioned in the introduction, Stoimenow and Hirasawa showed that the remaining 12--crossing
knots are fibered if and only if they satisfy Neuwirth's condition (\ref{abelianfib}). Altogether
this completes the classification of all fibered 12--crossing knots.

%Gabai \cite{Ga87} showed that 0-surgery on a knot in $S^3$ yields
%a fibered manifold if and only if the knot is fibered. We used
%Corollary \ref{corfibered} and \emph{KnotTwister} to show that
%0-surgery on any of the non--fibered knots with 12 crossings and
%monic Alexander polynomial is not fibered.

%======================================
\subsection{Closed manifolds} \label{sectionclosed} \label{sectionexamplesclosed}
Let $K\subset S^3$ be a non--trivial knot, denote the result of zero--framed surgery along $K$ by
$M_K$. Let $\phi \in H^1(M_K)$ be a generator. Gabai \cite[Theorem 8.8]{Ga87} showed that for a
non--trivial knot $K$ we have $||\phi||_{T,M_K}=2 \, \gen(K)-2$. Furthermore  Gabai \cite{Ga87}
showed that  a knot $K$ is fibered if and only if $M_K$ is fibered.

%Our calculations using  \emph{KnotTwister} show that the bounds from Theorem \ref{mainthm} for
%closed manifolds determine the Thurston norm for all $M_K$, where $K$ is a knot with 12 crossings
%or less. Furthermore our calculations show that twisted Alexander polynomials detect all
%non--fibered $M_K$, where $K$ is a knot with 12 crossings or less.

Using  \emph{KnotTwister} one can easily see that, for any knot $K$  with 12 crossings or less,
twisted Alexander polynomials corresponding to appropriate representations determine the Thurston
norm of $M_K$ and detect whether $M_K$ is fibered or not.
\section{Generalization of Cha's fibering obstruction}
\label{sectioncha}

In this section we formulate and prove Proposition
\ref{prop:rewrite} which, in combination with Theorem
\ref{mainthm2}, immediately implies Theorem \ref{corfibered2}. (In
Theorem \ref{corfibered2} it is easy to prove that if $\partial M$
is non--empty then $\partial M$ has to be  a collection of tori.)

In order to formulate Proposition \ref{prop:rewrite} we need some more notation. For a ring $R$ and
a maximal ideal $\m\subset R$ we denote the  field  $R/\m$ by $\F_\m$. Furthermore given a
representation $\a:\pi_1(M)\to \gl(R,k)$ we denote by $\a_\m$ the representation
$\pi_1(M)\xrightarrow{\a} \gl(R,k)\to \gl(\F_\m,k)$ where $\gl(R,k)\to \gl(\F_\m,k)$ is induced
from the canonical map $\pi_\m : R \to R/\m= \F_\m$. The main example to keep in mind is $R=\Z$,
$\m=(p)$ for a prime $p$, and $R/\m=\Z/(p) = \F_p$.

%Let  be a representation where $R$ is a Noetherian UFD. Recall that in this situation Cha
%\cite{Ch03} defined the twisted Alexander polynomial $\twialex \in R[t^{\pm 1}]$ which are
%well--defined up to multiplication by a unit in $R[t^{\pm 1}]$.

%Furthermore Cha showed that for a fibered knot $\twialex$ is monic. In this section we will
%generalize Cha's obstruction for fibering knots to obstructions for fibering (possibly
%\emph{closed}) 3--manifolds in Theorem \ref{corfibered2}.
%Cha's definition of
%twisted Alexander polynomials generalizes our definition.

\begin{proposition} \label{prop:rewrite}
Let $M$ be a 3--manifold whose boundary is empty or consists of tori and let $R$ be a Noetherian
UFD. Let $\phi\in H^1(M)$ be non--trivial and $\a:\pi_1(M)\to \gl(R,k)$ a representation. Then
$\twialex \in R\tpm$ is monic and
 \[ \tnphi
=\frac{1}{k}\degtaum.
%\big(\degtwialex-\degtwialexzero-\degtwialextwo \big).
\] if and only if for all maximal ideals $\m$ of $R$ we have that
$\Delta_1^{\a_\m}(t)$ is non--trivial and
\[ \tnphi
=\frac{1}{k}\deg(\tau(M,\phi,\a_\m))\]
%\left(\deg\left(\Delta_1^{\a_\m}(t)\right)-\deg\left(\Delta_0^{\a_\m}(t)\right)-\deg\left(\Delta_2^{\a_\m}(t)\right)
%\right) \]
\end{proposition}

\begin{proof}
%[Proof of Proposition \ref{prop:rewrite}]
We only prove this proposition in the case that $M$ is closed. The
proof for the case that $\partial M$ is a non--empty collection of
tori is very similar. Note that in either case $\chi(M)=0$.

We first make use of an argument in the proof of \cite[Theorem 5.1]{Mc02}. Choose a triangulation
$\tau$ of $M$. Let $T$ be a maximal tree in the $1$-skeleton of $\tau$ and let $T'$ be a maximal
tree in the dual 1-skeleton. We collapse $T$ to form a single 0-cell and join the 3-simplices along
$T'$ to form a single 3-cell. Denote the number of 1-cells by $n$. It follows from $M$ closed that
$\chi(M)=0$, hence there are $n$ 2-cells.  From the CW structure we obtain a chain complex
$C_*=C_*(\ti{M})$ of the following form
\[
0 \to C_3(\ti{M}) \xrightarrow{\partial_3} C_2(\ti{M})
\xrightarrow{\partial_2} C_1(\ti{M}) \xrightarrow{\partial_1}
C_0(\ti{M}) \to 0
\]
where $C_i(\ti{M})\cong \Z[\pi_1(M)]$ for $i=0,3$ and
$C_i(\ti{M})\cong \Z[\pi_1(M)]^n$ for $i=1,2$. Let $A_i,
i=0,\dots,3$ over $\Z[\pi_1(M)]$ be the matrices corresponding to
the boundary maps $\partial_i:C_i\to C_{i-1}$  with respect to the
bases given by the lifts of the cells of $M$ to $\ti{M}$. We can
arrange the lifts such that
 \[ \ba{rcl} A_3 &=&
(1-g_1, 1-g_2, \ldots, 1-g_n)^t,\\
A_1 &=& (1-h_1, 1-h_2, \ldots, 1-h_n). \ea \] Note that $\{g_1,\dots,g_n\}$ and $\{h_1,\dots,h_n\}$
are generating sets for $\pi_1(M)$ since $M$ is a closed 3--manifold. Since $\phi$ is non--trivial
there exist $r,s$ such that $\phi(g_r)\ne 0$ and $\phi(h_s)\ne 0$. Let $B_3$ be the $r$-th row of
$A_3$. Let $B_2$ be the result of deleting the $r$-th column and the $s$--th row from $A_2$. Let
$B_1$ be the $s$--th column of $A_1$.

Given a $p\times q$ matrix $B = (b_{rs})$ be  with entries in $\Z[\pi]$ we write $b_{rs}=\sum
b_{rs}^gg$ for $b_{rs}^g\in \Z, g\in \pi$. We then define $(\a\otimes \phi)(B)$ to be the $p\times
q$ matrix with entries $\sum b_{rs}^g \a(g)t^{\phi(g)}$. Since each $\sum b_{rs}^g
\a(g)t^{\phi(g)}$ is a $k\times k$ matrix with entries in $\ft$ we can think of  $(\a\otimes
\phi)(B)$ as a $pk\times qk$ matrix with entries in $\ft$.

Now note that
\[ \det((\a\otimes \phi)(B_3))=\det(\id-(\a\otimes \phi)(g_r))
=\det(\id-\phi(g_r)\a(g_r)) \ne 0 \]
 since $\phi(g_r)\ne 0$. Similarly $\det((\a\otimes
\phi)(B_1))\ne 0$ and $\det((\a_\m\otimes \phi)(B_i))\ne 0$, $i=1,3$ for any maximal ideal $\m$. We
need the following theorem which can be found in \cite{Tu01}.

\begin{theorem}\cite[Theorem~2.2, Lemma~2.5, Theorem~4.7]{Tu01} \label{thm:Tu22}
Let $S$ be a Noetherian UFD. Let $\b:\pi_1(M)\to \gl(S,k)$ be a
representation and $\varphi \in H^1(M)$. \bn \item If $
\det((\b\otimes \varphi)(B_i))\ne 0$ for $i=1,2,3$, then
$H_i^\b(M;S^k\tpm))$ is $S\tpm$--torsion for all $i$. \item If
$H_i^\b(M;S^k\tpm))$ is $S\tpm$--torsion for all $i$, and if
$\det((\b\otimes \varphi)(B_i))\ne 0$ for $i=1,3$, then
$\det((\b\otimes \varphi)(B_2))\ne 0$ and
\[ \prod\limits_{i=1}^3 \det((\b\otimes \varphi)(B_i))^{(-1)^{i}}
=\prod\limits_{i=0}^3 \left(\Delta_i^\b(t)\right)^{(-1)^{i+1}}
=\tau(M,\varphi,\b).\] \en
\end{theorem}

\noindent First assume that $\Delta_1^{\a_\m}(t)\ne 0$ and
\[ \tnphi
=\frac{1}{k}
\left(\deg\left(\Delta_1^{\a_\m}(t)\right)-\deg\left(\Delta_0^{\a_\m}(t)\right)-\deg\left(\Delta_2^{\a_\m}(t)\right)
\right)
\] for all maximal ideals $\m$. By Lemma \ref{lem:delta03} and
Proposition \ref{prop:delta2} we get $\Delta_i^{\a_\m}(t)\ne 0$ for all $i$, in particular
$H_i^{\a_\m}(M;\F_\m^k\tpm)$ is $\F_\m\tpm$--torsion for all $i$ and all maximal ideals $\m$. It
follows from Theorem \ref{thm:Tu22} that $\det((\a_\m\otimes \phi)(B_2))\ne 0$. Clearly this
also implies that $\det((\a \otimes \phi)(B_2))\ne 0$. Since we already know that $\det((\a
\otimes \phi)(B_i))\ne 0$ for $i=1,3$ it follows from Theorem \ref{thm:Tu22} that
$H_i^\a(M;R^k\tpm)$ is $R\tpm$--torsion for all $i$.

It follows from \cite[Lemma~4.11]{Tu01} that $\twialexzero$
divides $\det((\a\otimes \phi)(B_1))=\det(\id-\phi(h_s)\a(h_s))$
which is a monic polynomial in $R\tpm$ since $\phi(h_s)\ne 0$ and
since $\det(\a(h_s))$ is a unit.
%we have a representation $\a:\pi_1(M)\to \{ A\in \gl(R,k)| \det(A)\in R^*\}$.
But then $\twialexzero$ is monic as well. The same argument (again using \cite[Lemma~4.11]{Tu01})
shows that $\twialextwo$ is monic. It follows from the argument of Lemma \ref{lem:delta03} that
$H_3^\a(M;R^k\tpm)=0$, hence $\Delta_3^\a(t)=1$.

Denote the map $R \to R/\m=\F_\m$ by $\pi_\m$. We also denote the induced map $R\tpm \to \F_\m\tpm$
by $\pi_\m$. It follows from Theorem \ref{thm:Tu22} that
\[ \ba{rcl}
\prod\limits_{i=0}^3 \pi_\m\left(\Delta_i^\a(t)^{(-1)^{i+1}}\right) &=&
\prod\limits_{i=1}^3 \pi_\m\left(\det((\a\otimes \phi)(B_i))\right)^{(-1)^{i}}\\
&=&
\prod\limits_{i=1}^3 \det((\a_\m\otimes \phi)(B_i))^{(-1)^{i}}\\
&=&\prod\limits_{i=0}^3 \Delta_i^{\a_\m}(t)^{(-1)^{i+1}}\ea
\]
for all maximal ideals $\m$. By assumption we get
\[
\displaystyle\frac1k\sum\limits_{i=0}^3(-1)^{i+1}
\deg\big(\pi_\m\left(\Delta_i^\a(t)\right)\big)=\displaystyle\frac1k \sum\limits_{i=0}^3(-1)^{i+1}
\deg\left(\Delta_i^{\a_\m}(t)\right)= \tnphi \] for all $\m$. Since $\Delta_i^\a(t)$ is monic for
$i=0,2,3$ it follows that
\[
\deg\big(\pi_\m\left(\Delta_1^\a(t)\right)\big)=\deg\big(\pi_\n\left(\Delta_1^\a(t)\right)\big)\]
for all maximal ideals $\m$ and $\n$. Since $R$ is a UFD it follows that $\twialex$ is monic. Hence
$\deg\left(\pi_\m\left(\Delta_i^{\a}(t)\right)\right) = \deg\left(\Delta_i^{\a}(t)\right)$ for all
$i$ and all maximal ideals $\m$ and clearly
 \[ \tnphi
=\frac{1}{k} \big(\degtwialex-\degtwialexzero -\degtwialextwo
\big). \]

Now assume that $\twialex \in R\tpm$ is monic and
 \[ \tnphi
=\frac{1}{k} \big(\degtwialex-\degtwialexzero -\degtwialextwo
\big).
\] The same argument as above shows that $\Delta_i^\a(t)$, $i=0,2,3$, are monic as well.
Recall that $\det(\a\otimes \phi)(B_i)$, $i=1,3$, are monic
polynomials. It follows from Theorem \ref{thm:Tu22} that
\[ \det(\a\otimes \phi)(B_2)=
\det(\a\otimes \phi)(B_1)\, \det(\a\otimes \phi)(B_3)\,
\prod\limits_{i=0}^3 \left(\Delta_i^\a(t)\right)^{(-1)^{i+1}}\]
 is a quotient of monic non--zero polynomials. In particular
$\det(\a_\m\otimes \phi)(B_2)=\pi_\m(\det(\a\otimes \phi)(B_2))\ne 0$. It now follows immediately
from Theorem \ref{thm:Tu22} that $H_i^{\a_\m}(M;\F_\m^k\tpm))$ is $\F_\m\tpm$--torsion for all $i$.
In particular $\Delta_1^{\a_\m}(t) \ne 0$. Using arguments as above we now see that
\[ \ba{rcl}
\deg(\tau(M,\phi,\a_p))&=&\frac{1}{k}
\left(\deg\left(\Delta_1^{\a_\m}(t)\right)-\deg\left(\Delta_0^{\a_\m}(t)\right)
-\deg\left(\Delta_2^{\a_\m}(t)\right)
\right)\\
&=&\frac{1}{k}\sum\limits_{i=0}^3(-1)^{i+1}
\deg\left(\Delta_i^{\a_\m}(t)\right)\\
&=&\frac{1}{k}\sum\limits_{i=0}^3(-1)^{i+1}
\deg\left(\pi_\m\left(\Delta_i^\a(t)\right)\right)\\
&=&\frac{1}{k}\sum\limits_{i=0}^3(-1)^{i+1} \deg\left(\Delta_i^\a(t)\right)\\
&=&\tnphi. \ea \]
\end{proof}

\begin{remark}
Let $\a:\pi_1(M)\to \gl(\Z,k)$ be a representation. Then it is in general \emph{not} true that
for a prime $p$ we have $\Delta_1^{\a_p}(t)=\pi_p(\Delta_1^\a(t))\in \fpt$ (we use the notation
of Proposition \ref{prop:rewrite}), not even if $(M,\phi)$ fibers over $S^1$. Indeed, let $K$ be
the trefoil knot and $\varphi:\pi_1(X(K))\to S_3$ the unique epimorphism. Consider the
representation $\a(\varphi):\pi_1(X(K))\to \gl(\Z,2)$ as in Section \ref{sectionrep}. Then
$\deg\left(\pi_3(\Delta_1^\a(t))\right)=2$, but $\deg\left(\Delta_1^{\a_3}(t)\right)=3$.
\end{remark}

\end{document}